\documentclass[a4paper,11pt]{article}

\usepackage{amsfonts,amssymb,amsmath,amsthm,latexsym,epsf,epsfig,amscd,bbm,stmaryrd,stackrel}
\usepackage[english]{babel}
\usepackage{xcolor,colortbl}

\makeatletter \@addtoreset{equation}{section} \makeatother
\makeatletter \@addtoreset{enunciato}{section} \makeatother
\newcounter{enunciato}[section]

\newtheorem{ittheorem}{Theorem}
\newtheorem{itlemma}{Lemma}
\newtheorem{itproposition}{Proposition}
\newtheorem{itdefinition}{Definition}
\newtheorem{itremark}{Remark}
\newtheorem{itclaim}{Claim}
\newtheorem{itfact}{Fact}
\newtheorem{itconjecture}{Conjecture}
\newtheorem{itcorollary}{Corollary}

\newenvironment{theorem}{\addtocounter{enunciato}{1}
\begin{ittheorem}}{\end{ittheorem}}
\newenvironment{lemma}{\addtocounter{enunciato}{1}
\begin{itlemma}}{\end{itlemma}}
\newenvironment{proposition}{\addtocounter{enunciato}{1}
\begin{itproposition}}{\end{itproposition}}
\newenvironment{definition}{\addtocounter{enunciato}{1}
\begin{itdefinition}}{\end{itdefinition}}
\newenvironment{remark}{\addtocounter{enunciato}{1}
\begin{itremark}}{\end{itremark}}

\newenvironment{conjecture}{\addtocounter{enunciato}{1}
\begin{itconjecture}}{\end{itconjecture}}
\newenvironment{corollary}{\addtocounter{enunciato}{1}
\begin{itcorollary}}{\end{itcorollary}}
\newcommand{\be}[1]{\begin{equation}\label{#1}}
\newcommand{\ee}{\end{equation}}
\newcommand{\bl}[1]{\begin{lemma}\label{#1}}
\newcommand{\el}{\end{lemma}}
\newcommand{\br}[1]{\begin{remark}\label{#1}}
\newcommand{\er}{\end{remark}}
\newcommand{\bt}[1]{\begin{theorem}\label{#1}}
\newcommand{\et}{\end{theorem}}
\newcommand{\bd}[1]{\begin{definition}\label{#1}}
\newcommand{\ed}{\end{definition}}
\newcommand{\bp}[1]{\begin{proposition}\label{#1}}
\newcommand{\ep}{\end{proposition}}
\newcommand{\bc}[1]{\begin{corollary}\label{#1}}
\newcommand{\ec}{\end{corollary}}
\newcommand{\bcj}[1]{\begin{conjecture}\label{#1}}
\newcommand{\ecj}{\end{conjecture}}

\def \Z {{\mathbb Z}}
\def \R {{\mathbb R}}

\def \N {{\mathbb N}}

\def \ba {\begin{array}}
\def \ea {\end{array}}

\def \P {{\mathbb P}}
\def \E {{\mathbb E}}

\def \di {\mathrm{d}}

\def\1{\mathbbm{1}}

\begin{document}

\title{Continuity and anomalous fluctuations\\ in random walks in
\\dynamic random environments:
\\numerics, phase diagrams and conjectures.}

\author{\renewcommand{\thefootnote}{\arabic{footnote}}
L.\ Avena \footnotemark[1]
\\
\renewcommand{\thefootnote}{\arabic{footnote}}
P.\ Thomann \footnotemark[1]
}

\footnotetext[1]{Institut f\"ur Mathematik, Universit\"at
Z\"urich, Winterthurerstrasse 190, Z\"urich, CH- 8057,
Switzerland. \\ E-mail: luca.avena@math.uzh.ch.}

\maketitle

\begin{abstract}
We perform simulations for one dimensional continuous-time random
walks in two dynamic random environments with fast (independent
spin-flips) and slow (simple symmetric exclusion) decay of
space-time correlations, respectively. We focus on the asymptotic
speeds and the scaling limits of such random walks. We observe
different behaviors depending on the dynamics of the underlying
random environment and the ratio between the jump rate of the
random walk and the one of the environment. We compare our data
with well known results for static random environment. We observe
that the non-diffusive regime known so far only for the static
case can occur in the dynamical setup too. Such anomalous
fluctuations give rise to a new phase diagram. Further we discuss
possible consequences for more general static and dynamic random
environments.

\vspace{0.5cm}\noindent

\end{abstract}

\newpage


\section{Introduction}
\label{S1}


\subsection{Random walk in static and dynamic random environments}
\label{S1.1}

Random Walks in Random Environments (RWRE) on the integer lattice
are RWs on $\Z^d$ evolving according to random transition kernels,
i.e., their transition probabilities/rates depend on a random
field (\emph{static} case) or a stochastic process (\emph{dynamic}
case) called Random Environment (RE). Such models play a central
role in the field of disordered systems of particles. The idea is
to model the motion of a particle in an inhomogeneous medium. In
contrast with standard homogeneous RW, RWRE may show several
unusual phenomena as non-ballistic transience, non-diffusive
scalings, sub-exponential decay for large deviation probabilities.
All these features are due to impurities in the medium that
produce trapping effects. Although they have been intensively
studied by the physics and mathematics communities since the early
70's, except for the one-dimensional static case and few other
particular situations, most of the results are of qualitative
nature, and their behavior is far from being completely
understood. We refer the reader to \cite{Sz02,Ze06} and
\cite{Av10,DoKeLi08} for recent overviews of the state of the art
in \emph{static} and \emph{dynamic} REs, respectively.

In this paper we focus on two one-dimensional models in dynamic
RE. In particular, the RW will evolve in continuous time in
(two-states) REs given by two well known interacting particle
systems: independent spin flip and simple symmetric exclusion
dynamics. Several classical questions regarding these types of
dynamical models are still open while the behavior of the
analogous RW in a i.i.d.\ static case is completely understood. We
perform simulations focusing on their long term behavior. We see
how such asymptotics are influenced as a function of the jump rate
$\gamma$ of the dynamic REs. The idea is that by tuning the speed
of the REs, we get close to the static RE ($\gamma$ close to $0$)
or to the averaged medium ($\gamma$ approaching infinity). We
observe different surprising phases which allow us to set some new
challenging conjectures and open problems.

Although the choice of the models could appear too restrictive and
of limited interest, as it will be clear in the sequel, such
particular examples present all the main features and the rich
behavior of the general models usually considered in the RWRE
literature. The conjectures we state can be extended to a more
general setup (see Section \ref{S4.3}).

The paper is organized as follows. In this section we define the
models and give some motivation. In Section \ref{S2} we review the
results known for the analogous RW in a static RE. Section
\ref{S3} represents the main novel. We present therein the results
of our simulations which shine a light on the behavior of the
asymptotic speed (Section \ref{S3.1}) and on the scaling limits of
such processes (Sections \ref{S3.2} and \ref{S3.3}). When
discussing each question we list several conjectures. In the last
Section \ref{S4} we present a brief description of the algorithms,
we discuss the robustness of our numerics and possible
consequences for more general RE.


\subsection{The model}
\label{S1.2} We consider a one-dimensional RW whose transition
rates depend on a dynamic RE given by a \emph{particle system}. In
Section \ref{S1.2.1}, we first give a rather general definition of
\emph{particle systems} and then we introduce the two explicit
examples we will focus on. In Section \ref{S1.2.2} we define the
RW in such dynamic REs.

\subsubsection{Random Environment: particle systems}
\label{S1.2.1} Let $\Omega=\{0,1\}^{\Z}$. Denote by
$D_\Omega[0,\infty)$ the set of paths in $\Omega$ that are
right-continuous and have left limits. Let $\{P^\eta,
\eta\in\Omega\}$ be a collection of probability measures on
$D_\Omega[0,\infty)$. A
\emph{particle system} \be{xidef} \xi = (\xi_t)_{t \geq 0} \quad
\mbox{ with } \quad \xi_t = \{\xi_t(x)\colon\,x\in\Z\}, \ee is a
Markov process on $\Omega$ with law $P^\eta$, when
$\xi_0=\eta\in\Omega$ is the starting configuration. Given a
probability measure $\mu$ on $\Omega$, we denote by $P^\mu(\cdot)
:= \int_{\Omega} P^\eta(\cdot)\,\mu(\di\eta)$ the law of $\xi$
when $\xi_0$ is drawn from $\mu$. We say that site $x$ is occupied
by a particle (resp. vacant) at time $t$ when $\xi_t(x)=1$ (resp.
$0$).

Informally, a \emph{particle system} is a collection of particles ($1$'s)
on the integer lattice evolving in a Markovian way. Depending on
the specific transition rates between the different
configurations, one obtains several types of \emph{particle
systems}. Each particle may interact with the others: the
evolution of each particle is defined in terms of local transition
rates that may depend on the state of the system in a neighborhood
of the particle. For a formal construction, we refer the reader to
Liggett \cite{Li85}, Chapter I.
\bigskip

In the sequel we focus on two well known examples with strong and
weak mixing properties, respectively.

\bigskip
(1) {\bf{Independent Spin Flip (ISF)}}

\bigskip
Let $\xi = (\xi_t)_{t\geq 0}$ be a one-dimensional
\emph{independent spin-flip system}, i.e., a Markov process on
state space $\Omega$ with generator $L_{ISF}$ given by
\be{Generator} (L_{ISF}f)(\eta) = \sum_{x\in\Z}
\left[\lambda\eta(x)+\gamma\left(1-\eta(x)\right)\right][f(\eta^x)-f(\eta)],
\qquad \eta\in\Omega, \ee where $\lambda,\gamma\geq0$, $f$ is any
cylinder function on $\Omega$, $\eta^x$ is the configuration
obtained from $\eta$ by flipping the state at site $x$.

In words, this process is an example of a \emph{non-interacting
particle system} on $\{0,1\}^{\Z}$ where the coordinates
$\eta_t(x)$ are independent two-state Markov chains, namely, at
each site (independently with respect to the other coordinates)
particles flip into holes at rate $\lambda$ and holes into
particles at rate $\gamma$. This particle system has a unique
ergodic measure given by the Bernoulli product measure with
density $\rho=\gamma/(\gamma+\lambda)$ which we denote by
$\nu_\rho$ (see e.g. \cite{Li85}, Chapter IV).

\bigskip

(2) {\bf{Simple Symmetric Exclusion (SSE)}}

The SSE is an \emph{interacting particle system} $\xi$ in which
particles perform a simple symmetric random walk at a certain rate
$\gamma>0$ with the restriction that only jumps on vacant sites
are allowed. Formally, its generator $L_{SSE}$, acting on cylinder
functions $f$, is given by \be{SSEgen} (L_{SSE}f)(\eta) =
\gamma\sum_{ {x,y\in\Z} \atop {x \sim y} }
[f(\eta^{x,y})-f(\eta)], \qquad \eta\in\Omega, \ee where  the sum
runs over unordered neighboring pairs of sites in $\Z$, and
$\eta^{x,y}$ is the configuration obtained from $\eta$ by
interchanging the states at sites $x$ and $y$.

It is known (see \cite{Li85}, Chapter VIII) that the family of
Bernoulli product measures $\nu_\rho$, with density $\rho\in(0,1)$,
characterize the set of equilibrium measures for this dynamics.

\bigskip

\br{Remark}\label{remark} Note that the ISF and the SSE are
completely different types of dynamics. They are both Markovian in
time but while the ISF has no spatial correlations, the SSE has
space-time correlations. The ISF has very good mixing properties
due to the fact that once $\gamma$ and $\lambda$ are given, no
matter of the starting configuration, it will converge
exponentially fast to the unique equilibrium given by $\nu_\rho$
with $\rho=\gamma/(\gamma+\lambda)$. On the contrary, the SSE
dynamics is strongly dependent on the starting configuration and
therefore does not satisfy any uniform mixing condition. In fact,
it is a conservative type of dynamics with a family of equilibria
given by $\left\{\nu_\rho : \rho\in(0,1)\right\}$. Because of
these substantial differences, in the sequel we will informally say
that the ISF and the SSE are examples of \emph{fast} and
\emph{slowly mixing} dynamics, respectively. \er

\subsubsection{RW on particle systems}
\label{S1.2.2}

Conditional on the particle system $\xi$, let \be{rwdef} X =
(X_t)_{t\geq 0} \ee be the continuous time random walk jumping at
rate $1$ with local transition probabilities \be{rwtrans}
\begin{aligned}
&x \to x+1 \quad \mbox{ at rate } \quad p\,\xi_t(x) + (1-p)\,[1-\xi_t(x)],\\
&x \to x-1 \quad \mbox{ at rate } \quad (1-p)\,\xi_t(x) + p\,[1-\xi_t(x)],
\end{aligned}
\ee
with
\be{p}
p\in[1/2,1).\ee

In words, the RW $X$ jumps according to an exponential clock with
rate $1$, if $X$ is on occupied sites (i.e. $\xi_t(X_t)=1$), it
goes to the right with probability $p$ and to the left with
probability $1-p$, while at vacant sites it does the opposite.

We write $P^\xi_0$ to denote the law of $X$ starting from $X(0)=0$
conditional on $\xi$, and \be{Pxnudef} \P_{\mu,0}(\cdot) =
\int_{D_\Omega[0,\infty)} P^\xi_0(\cdot)\,P^\mu(\di\xi) \ee to
denote the law of $X$ averaged over $\xi$. We refer to $P^\xi_0$
as the \emph{quenched} law and to $\P_{\mu,0}$ as the
\emph{annealed} law. In what follows, when needed, we will denote
by \be{Xparameters}X(p,\gamma,\rho),\ee the RW $X$ just defined
either in the ISF or in the SSE environment starting from
$\nu_\rho$ and jumping at rate $\gamma$. Note that in the ISF
case, the parameter $\lambda$ is uniquely determined once we fix
$\gamma$ and $\rho$.

From now on we assume w.l.o.g. $\rho\in[1/2,1)$. The choice of
$p,\rho\in[1/2,1)$ is not restrictive, indeed, due to symmetry, it
is easy to see the following equalities in distribution
\be{symmetry}X(p,\rho,\gamma)\stackrel{\P}{=}X(1-p,1-\rho,\gamma)
\stackrel{\P}{=}-X(p,1-\rho,\gamma).\ee


\subsection{On mixing dynamics}
\label{S1.3} In our models, the REs at each site have only two
possible states ($0$ or $1$), in most of the literature on RWRE,
the models are defined in a more general framework where
infinitely many states are allowed. The first paper dealing with
RW in dynamic RE goes back to \cite{BoIgMaPe92}. Since then, there
has been intensive activity and several advances have recently
been made showing mostly LLN, invariance principles and LDP under
different assumptions on the REs or on the transition
probabilities of the walker. See for example
\cite{AvdHoRe11,BaZe06,BoMiPe00,BoMiPe04,BrKu09,DoKeLi08,RaAgSe05}
(most of these references are in a discrete-time setting). For an
extensive list of reference we refer the reader to
\cite{Av10,DoKeLi08}.

One of the main difficulties in the analysis of random media arises
when space-time correlations in the RE are allowed. Both models
presented in Section \ref{S1.2.1} fit in this class but, as
mentioned in Remark \ref{remark}, their mixing properties are
substantially different. The ISF dynamics belongs to the class of
\emph{fast mixing} environments which is known to be qualitatively
similar to a homogeneous environment. In fact, a RW on this type
of RE exhibits always diffusive scaling.

The SSE is an example of what we called \emph{slowly mixing}
dynamics. For a RW driven by these latter types of dynamics, we
are not aware of any results other than \cite{AvdHoRe10,AvSaVo11,dHoKeSi12}.
One of the main result of our simulations is that the RW in
\eqref{rwdef} on the SSE, similarly to the RW in a static RE (see
Section \ref{S2.2}), may exhibit non-diffusive behavior (see
Section \ref{S3.2}). This latter result is related to trapping
phenomena (see Sections \ref{S2.4} and \ref{S2.6}) and suggests
that, when considering \emph{non-uniform slowly mixing}
environments, the medium looks substantially different than a
homogeneous one. This is confirmed by the rigorous annealed large
deviation results in \cite{AvdHoRe10} (see the last paragraph in
Section \ref{S2.3}). Note that these results are due to the
correlation structure of the RE but also depend on the following
essential ingredients which produce some strong trapping effect:
the one dimensional setting, the RW and the SSE being both
nearest-neighbor, and the presence of local drifts to the right
and to the left for the RW.

\subsection{Particle systems as random environments}
\label{S1.4} The reader may wonder why we consider random
environments given by particle systems. Particle systems represent
a natural physical example of a two-state dynamical RE. Particle
system theory has been intensively developed in the last thirty
years and results from this theory can be used in the present
context (see e.g.
\cite{Av10,AvdHoRe10,AvdHoRe11,AvSaVo11,dHoKeSi12,dHoSaSi11,ReVo11}).
Several results proven for such particular models can be extended
to more general settings. Finally, as shown in this paper, these
dynamics are not too complex from an algorithmic point of view
allowing to obtain good approximations for the asymptotics.


\section{Static case and trapping phenomena: review }
\label{S2} We present in this section well known results for the
analogous model in an i.i.d.\ static medium. Consider a
\emph{static} random environment $\eta\in\{0,1\}^\Z$ with law
$\nu_\rho$, the Bernoulli product measure with density $\rho \in
[1/2,1)$. Given a realization of $\eta$, let $X=(X_t)_{t \geq 0}$
be the random walk with transition rates (compare with
(\ref{rwtrans})) \be{rwtransdis}
\begin{aligned}
&x \to x+1 \quad \mbox{ at rate } \quad p\eta(x)+(1-p)[1-\eta(x)]=c^+(\eta),\\
&x \to x-1 \quad \mbox{ at rate } \quad (1-p)\eta(x)+p[1-\eta(x)]=c^-(\eta),
\end{aligned}
\ee
with
\be{p2}
p\in[1/2,1).\ee

\subsection{Recurrence and LLN }
\label{S2.1}

\begin{figure}[hbtp]
\begin{center}
\includegraphics[width=.6\textwidth]{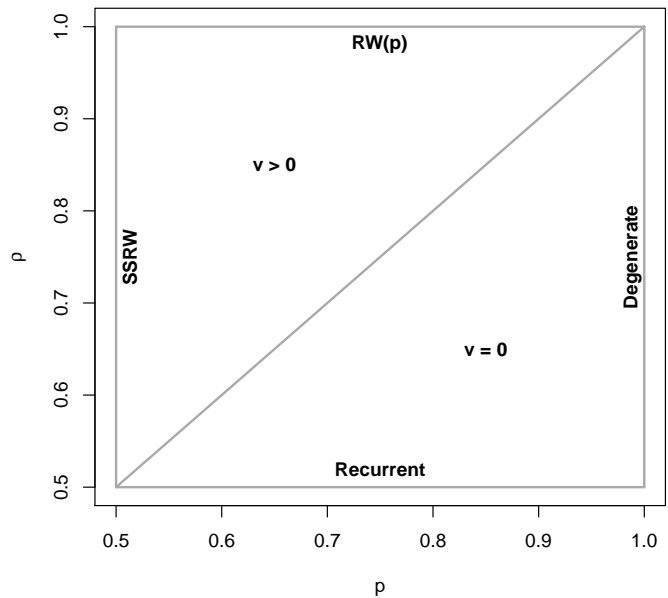}
\end{center}
\caption{ The sides of this square represent degenerate cases. In
particular when $p=1/2$ or $\rho=1$, the RW $X$ does not feel
anymore the environment behaving as a Simple Symmetric Random Walk
(SSRW in the picture) or as a homogeneous RW with drift $2p-1$
(RW(p) in the picture), respectively. When $p=1$ we are in a
trivial degenerate case. When $\rho=1/2$ we are in the
recurrent case. Inside the square, by \eqref{SignStSpeed}, above
the diagonal we have transience with positive speed, while at and
below the diagonal a non-ballistic transient regime holds, i.e.\ transience
at zero-speed.} \label{StaticSpeed}
\end{figure}

In \cite{So75} it is shown that $X$ is recurrent when $\rho
=\tfrac12$ and transient to the right when $\rho>\tfrac12$. In the
transient case both ballistic and non-ballistic behavior occur
(see Figure \ref{StaticSpeed}), namely, $\lim_{t\to\infty}
X_t/t=v_{\mathrm{static}}$ exists for $\P_{\nu_\rho}$-a.e.\
$\eta$, and

\be{SignStSpeed}
v_{\mathrm{static}}\,\,\left\{\begin{array}{ll}
= 0 &\mbox{ if } \rho \in [\frac12,p],\\
> 0 &\mbox{ if } \rho \in (p,1].
\end{array}
\right.\ee

In particular, for $\rho \in (p,1]$, \be{StSpeed1} v_{\mathrm{static}} =
v_{\mathrm{static}}(\rho,p) =
(2p-1)\,\frac{\rho-p}{\rho(1-p)+p(1-\rho)}. \ee

\newpage
\subsection{Scaling limits}
\label{S2.2}

The scaling limits of one-dimensional RWRE have been derived in
quite a general framework in a series of papers (see
\cite{KeKoSp75,Si82} and \cite{Ze04,Ze06} for a review of those
results and references). It turns out that diffusive,
super-diffusive or sub-diffusive regimes can occur. By diffusive
regime, we mean that $X_t-vt$ divided by $\sqrt{t}$ converges in
distribution to a non-degenerate Gaussian distribution, while we
refer to super- or sub-diffusive regime when $X_t-vt$ has to be
rescaled by some factor of order $t^\alpha$, with $\alpha$ bigger
or smaller than $1/2$, respectively, to converge weakly to some
non-degenerate distribution. We now review the different scalings
in the i.i.d. static RE. In what follows we focus only on the
annealed law.

When $X$ is recurrent, \cite{Si82} showed that $X$ is extremely
sub-diffusive and it converges weakly to a non-degenerate random
variable $Z$, namely, \be{Sinai}\frac{\sigma^2 X_t}{(\log
t)^2}\stackrel[t\to\infty]{(\P_{\nu_\rho})}{\longrightarrow}Z, \ee
where $\sigma^2$ is a positive constant and $Z$ is a random
variable with a non-trivial law that was later identified by
Kesten \cite{Ke86}. In this case, $X$ is called Sinai's random
walk.

\begin{figure}[hbtp]
\begin{center}
\includegraphics[width=.7\textwidth]{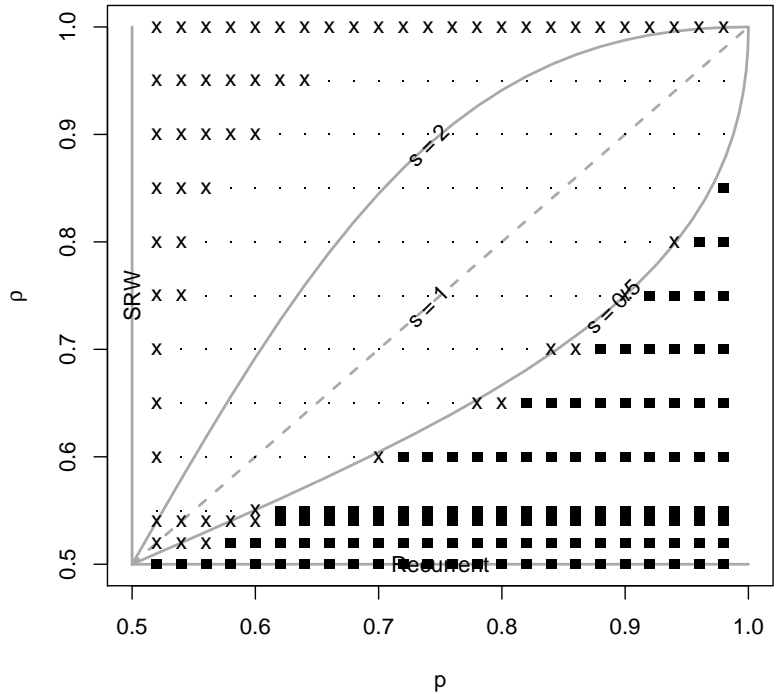}
\end{center}
\caption{For each of the points marked with a black square, a dot
or a cross, we performed numerical experiments (similar to those
explained in Section \ref{S3.3}) to determine the scaling
exponents of $X=X(p,\rho)$. The symbols square, dot or cross mean
that for the corresponding points, our numerical estimates gave a
sub-, a super- or a diffusive scaling exponent, respectively. Note
that except for a small region in between super-diffusive and
diffusive regimes (it is reasonable to have numerical fluctuations
close to a phase transition), the experiments confirm the
theoretical picture. Figure \ref{ScalingDataStatic2} provides a
few explicit examples of our numerics in the different regimes.}
\label{ScalingDataStatic1}
\end{figure}

\begin{figure}[hbtp]
\begin{center}
\includegraphics[width=.8\textwidth]{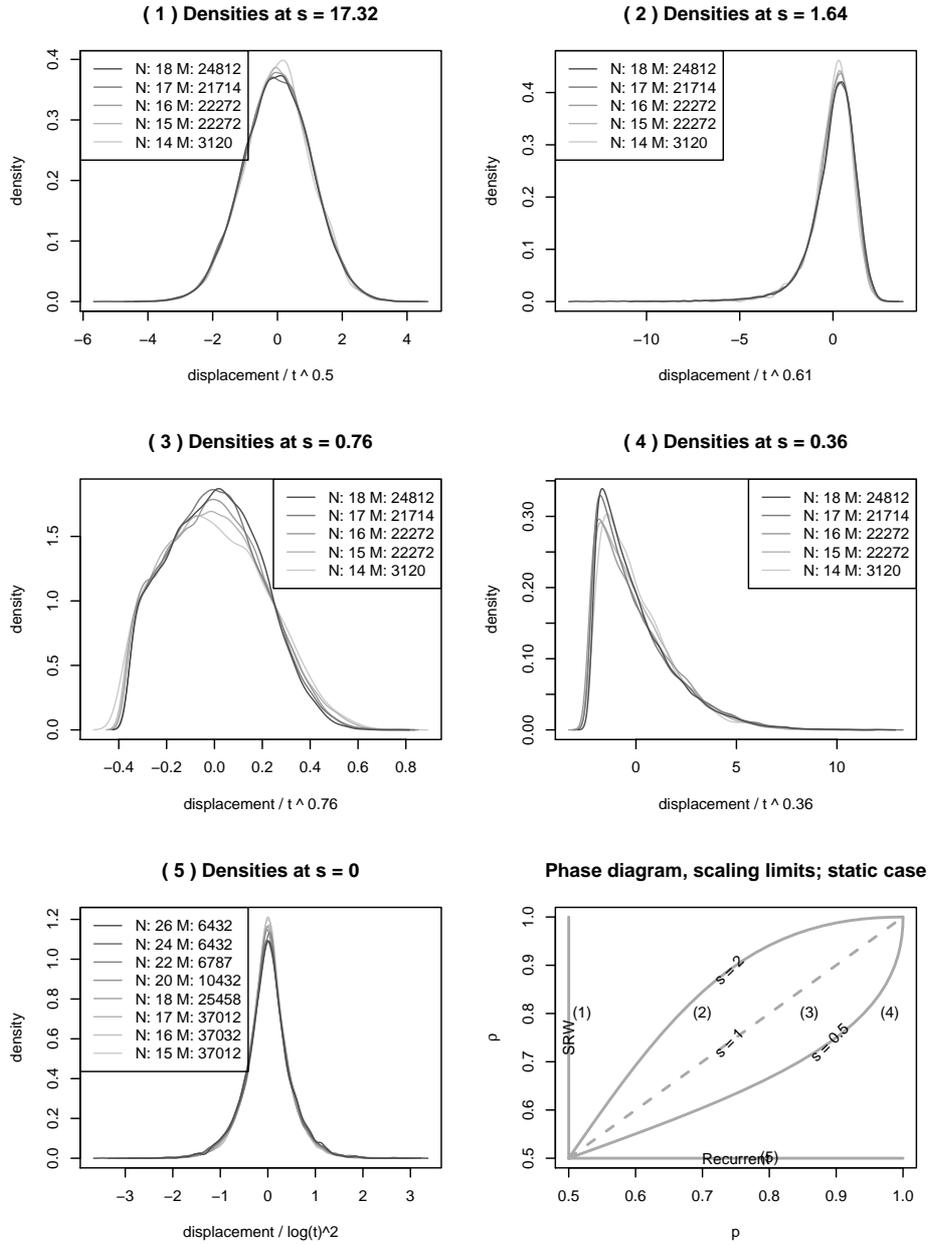}
\end{center}
\caption{A few explicit examples of the densities obtained via our
numerics in each of the different scaling regimes. In the right
bottom picture the $(p,\rho)$-points associated to each labelled
plot are specified. In each plot we overlapped the densities
obtained with independent experiments at the different times
$n=2^N$, over samples of size $M$, properly rescaled. In
particular, (1),(2),(3),(4) and (5) correspond to ballistic
diffusive, ballistic super-diffusive, transient zero-speed
super-diffusive, transient zero-speed sub-diffusive and recurrent
Sinai case, respectively.} \label{ScalingDataStatic2}
\end{figure}

When $X$ is right-transient, \cite{KeKoSp75} proved that the key
quantity to determine the right scaling is the root $s$ of the
equation
\be{SSS}\E_{\nu_{\rho}}\left[\left(\frac{c^-(\eta)}{c^+(\eta)}
\right)^s\right]=1,\ee where $c^-(\eta)$ and $c^+(\eta)$ represent
the rates to jump to the left and to the right, respectively (see
\eqref{rwtransdis}). In particular, they proved that when $s>2$,
$X$ is diffusive with Gaussian limiting distribution, while for
$s\in(0,2]$ super- or sub-diffusivity occur with some non-trivial
stable law of parameters $(s,b)$ as limiting distribution ($b$ is
a constant, see also Theorem 2.3 in \cite{Ze06} for more details
and references). The proof is based on the analysis of hitting
times and makes use of the extra assumption that
$\log\left[\frac{c^-(\eta)}{c^+(\eta)}\right]$ has a
non-arithmetic distribution. This latter is a delicate technical
assumption (see also Remark 3 in \cite{KeKoSp75}) not satisfied in
our model since
\be{arithmetic}\log\left[\frac{c^-(\eta)}{c^+(\eta)}\right]=
[2\eta(0)-1]\log\left(\frac{1-p}{p}\right)\ee does have an
arithmetic distribution. At the present state of the art, the role
of this assumption is not entirely clear. There are examples in
the literature in which by dropping it, the convergence in
distribution does not hold (see e.g.\ \cite{Bo06} Section 8 and
references therein). For our model \eqref{rwtransdis}, we
performed simulations (see Figures \ref{ScalingDataStatic1} and
\ref{ScalingDataStatic2}) which clearly suggest that the
arithmetic law of \eqref{arithmetic} does not play any role,
namely, the scaling behavior of $X$ is like in the general case
under the assumption of a non-arithmetic law. The following list
summarizes the different scalings of $X$.

\begin{itemize}
\item Diffusive: $s>2$, scaling order $\sqrt{t}$.
\item Super-diffusive: $s\in(0.5,2)$.
\begin{itemize}
\item $s\in(1,2)$: scaling order $t^{1/s}$.
\item $s=1$ :      scaling order $t/\log t$.
\item $s\in(0,1)$:  scaling order $t^s$.
\end{itemize}
\item Diffusive: $s=0.5$, scaling order $\sqrt{t}$.
\item Sub-diffusive:   $s\in(0,0.5)$, scaling order $t^s$.
\end{itemize}

In our model, the explicit solution of \eqref{SSS} is given by
\be{Sscaling}s=s(p,\rho)=
\frac{\log\left(\frac{1-\rho}{\rho}\right)}{\log\left(\frac{1-p}
{p}\right)}>0, \quad \text{for } p,\rho>1/2.
\ee

\begin{figure}[hbtp]
\begin{center}
\includegraphics[width=.75\textwidth]{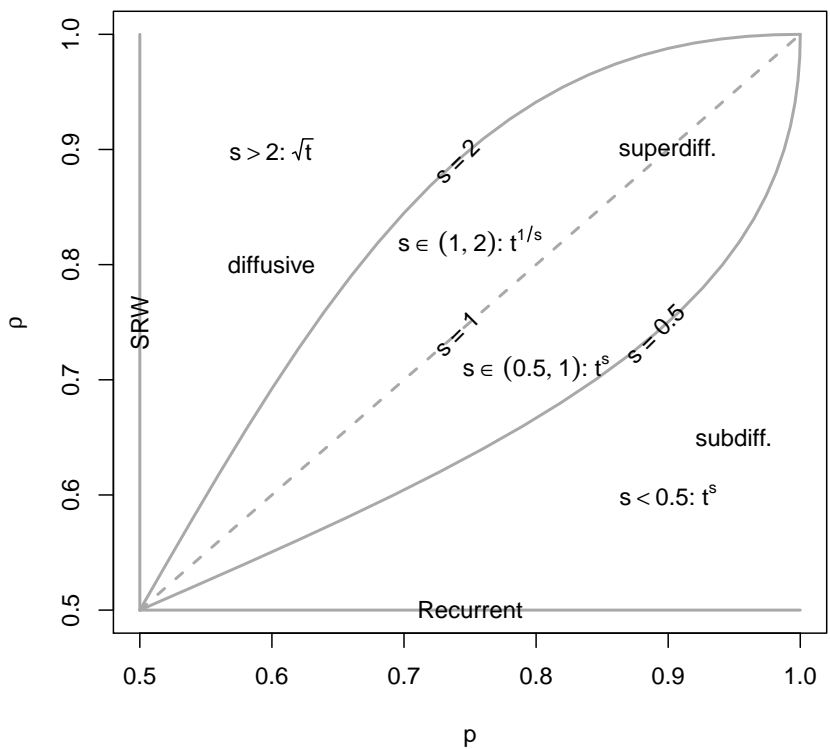}
\end{center}
\caption{ The sides of the square represent degenerate cases (see
Figure \ref{StaticSpeed}). When $\rho=1/2$ we are in the so called
Sinai case in which $X$ is recurrent and extremely sub-diffusive.
We call ``leaf'' the region around the diagonal delimited by the
curves $f_1(p)=(1-p)^2/p^2$ and $f_2(p)=(p-\sqrt{p(1-p)})/(2p-1)$
corresponding to $s=2$ and $s=0.5$, respectively. The area above
the leaf corresponds to the diffusive case with $s>2$, inside the
leaf we have the super-diffusive regime while in the lower
remaining region for $s\in(0,0.5)$ we have the sub-diffusive
regime.} \label{ScalingStatic}
\end{figure}

Figure \ref{ScalingStatic} shows the phase diagram in $(\rho,p)$
of the regimes just described.

\newpage

As we mentioned, we tested numerically the results presented so
far (see Section \ref{S4.1} for a description of the algorithms we
implemented). Figures \ref{ScalingDataStatic1} and
\ref{ScalingDataStatic2} show that our numerics match the
theoretical picture just described.

\newpage

\subsection{Large Deviations for the empirical speed}
\label{S2.3} For the empirical speed of a one-dimensional RW in a
static RE, quenched and annealed Large Deviation Principles (LDPs)
and refined large deviations estimates have been obtained in a
series of papers (see e.g. \cite{CoGaZe00,GedHo94,Va03,Ze06}). We
just mention that for the RW $X$ defined through
\eqref{rwtransdis}, when $v>0$, the rate function associated to
the LDP at rate $t$ is zero on the whole interval $[0,v]$. Roughly
speaking, this is saying that for $\theta\in[0,v)$,
$P(X_t/t\approx \theta)$ decays sub-exponentially in $t$. We recall
that for homogeneous RW such a decay is always exponential. This
phenomenon is due to trapping effects which we briefly introduce
in the next section.

Large deviation estimates for the dynamical models in Section
\ref{S1.2} were obtained in \cite{AvdHoRe10}. In particular, it is
shown that in the ISF case the rate function has a unique zero (as
for homogeneous RW) while in the SSE case, the rate function (at
least under the annealed measure) presents a flat piece as we just
described for the i.i.d. static case.

\subsection{Trapping phenomena}
\label{S2.4} The anomalous behaviors like the transient regime at
zero-speed, the non-diffusivity, as well as the sub-exponential
decay of the large deviations probabilities we reviewed, are due
to the presence of traps in the medium, i.e., localized regions
where the walk spends a long time with a high probability. To get
an intuition, the next picture shows an explicit example of a
trap. For a deeper insight of trapping phenomena we should
introduce the random potential representation of the environment
for which we refer the reader to the literature (see e.g.
\cite{Si82} and other references in \cite{Ze06}).

\begin{figure}[hbtp]
\begin{center}
\includegraphics[width=8cm]{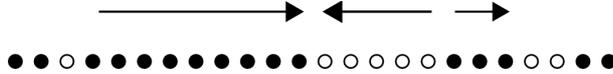}
\end{center}
\caption{An example of a trap, even though the global drift is to
the right, i.e. $\rho>1/2$, a long interval with vacant sites
creates a region with local drift against the global one. To cross
such a trap $X$ needs an average number of trials that is
exponential in the size of the interval.} \label{Trap}
\end{figure}

\subsection{The averaged medium}
\label{S2.5}
In the sequel we will refer to the RW in the \emph{averaged medium},
$Y(\mathrm{averaged})$, for the homogeneous nearest neighbor RW
with transition rates
\be{AvMeRW}
\begin{aligned}
&x \to x+1 \quad \mbox{ at rate } \quad p\,\rho + (1-p)\,(1-\rho),\\
&x \to x-1 \quad \mbox{ at rate } \quad (1-p)\,\rho + p\,(1-\rho),
\end{aligned}
\ee where $\rho,p\in[1/2,1)$. It easily follows that such a RW is
right transient as soon as $p>1/2$ and $\rho\neq1/2$. Moreover, by
the law of large numbers for i.i.d. sequences, we get that
\be{AvMedSpeed}\begin{aligned}\lim_{t\to\infty}
\frac{Y_t(\mathrm{averaged})}{t}&=v_{\mathrm{averaged}}(p,\rho)\\
&=(2\rho-1)(2p-1),\;a.s.
\end{aligned}
\ee

This RW would correspond to the walker defined in \eqref{rwdef}
which observes at each lattice position a constant density of
particles $\rho$. That is, between two jumps of the walker the
environment is replaced by one that is an independent sample drawn
from $\nu_\rho$.

In the dynamical models of Section \ref{S1.2.1}, both particle
systems are assumed to be in equilibrium with density $\rho$ and
exhibit some decay of correlations in time. Hence, roughly
speaking, if $\gamma$ approaches infinity, the environment becomes
asymptotically independent between two jumps of the walker.
Therefore, under the annealed law, we expect for both models that 
there is some convergence to the averaged medium as $\gamma \to\infty$.

\subsection{Towards the dynamic RE: dissolvence of traps}
\label{S2.6}

In the previous sections we saw that the RW $X$ in the static RE
$\eta\in\{0,1\}^\Z$ sampled from the Bernoulli product measure
$\nu_\rho$ presents ``slow-down phenomena" due to the presence of
\emph{traps}. The dynamical models in Section \ref{S1.2.2} can be
interpreted as the model in the static RE when we ``switch on''
some stochastic dynamics which allows particles of the RE to move
(SSE dynamics) or to be created/annihilated (ISF dynamics). The
natural question is then:
$$\text{How does the dynamics of the random environment influence}$$
$$\text{the trapping effects present in the static case?}$$

In the sequel we will present the outcome of simulations for the
empirical speed of $X$ in the different REs of Section
\ref{S1.2.1} and we will compare them with the static and the
averaged medium case. Note that  $\nu_\rho$ is an equilibrium
measure for both particle systems we use: ISF and SSE.

\begin{figure}[hbtp]
\begin{center}
\includegraphics[width=11cm]{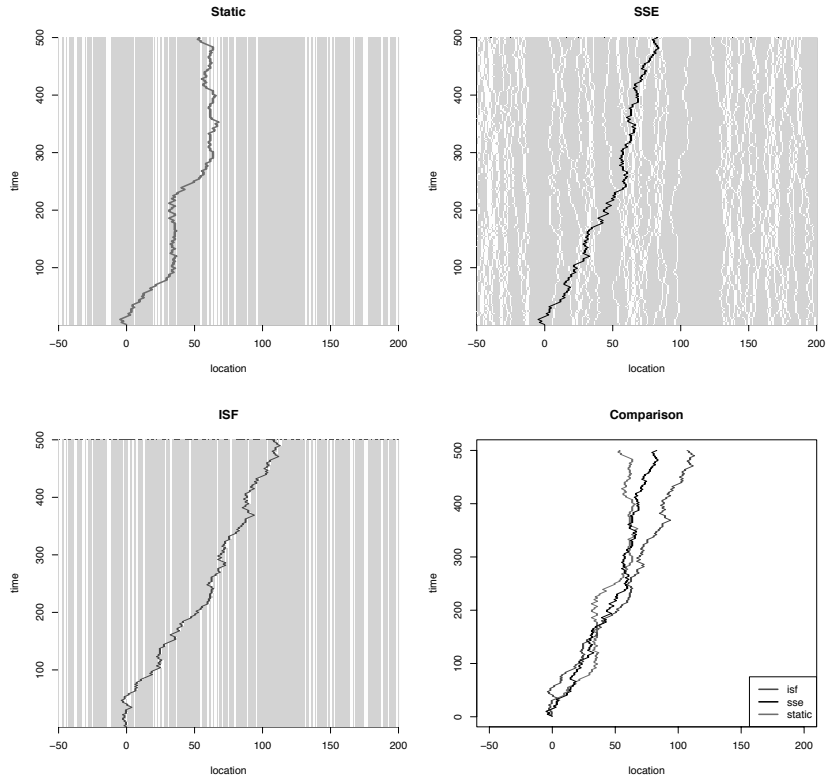}
\end{center}
\caption{Three simulated trajectories of $X$ in the static, ISF
and SSE environments starting from the same configuration.
$(p,\rho)=(0.7,0.8)$, and $\gamma=0.1$ for ISF and SSE. These
trajectories are compared in the bottom right picture showing the
intuitive domination mentioned above. The background of the other
plots displays realizations of the corresponding RE. In particular
gray and white mean presence or absence of particles,
respectively. Note that in the ISF case, due to the independence
in space of the dynamics, the RE has been updated only around the
RW trajectory, (see the description of the algorithm in Section
\ref{S4.1}).} \label{Trajectories}
\end{figure}

At a heuristic level, one should expect that the evolution of
particles in dynamic REs favors the dissolvence of traps,
consequently the RW $X$ in the dynamic RE should be ``faster''
than in the static medium. In other words, the long stretches of
holes present at time zero, and responsible of the slow-down
phenomena in the static case, get destroyed by the movement of
particles in the dynamic case. As a counter effect the dynamics
can create new traps during the evolution of the RW $X$.
Nevertheless, in the static case, the \emph{traps} are frozen,
while in the dynamic case, all the \emph{traps} have an a.s.\
finite survival time. Such intuitive arguments suggest that the
displacement of the RW $X$ should be bigger in the dynamic case
than in the static one.
Figure \ref{Trajectories} illustrates this intuitive domination; it
represents simulated trajectories of the RW in the three different
random environments (static, ISF, SSE) starting from the same
configuration sampled from $\nu_\rho$ at a given $p$.

Furthermore, depending on the specific dynamics of the underlying
particle system, the survival time and the nature of a typical trap
have to be different.
In fact, if we consider a trap as in Figure \ref{Trap} formed by
an interval filled of holes. It is clear that in the ISF case,
since particles can be created at each site at a given rate,
such a trap gets easily destroyed. In the SSE cases, due to the
conservation law, to dissolve such a trap, we have to wait for
particles from outside the interval to invade it.

\newpage


\section{Results and conjectures}
\label{S3} In this section we finally present the outcome of the
simulations for the asymptotics of $X$ in the different cases.
Section \ref{S3.1} concerns the asymptotic speed as a function of
the parameters $(\rho,p,\gamma)$ while Sections \ref{S3.2} and
\ref{S3.3} focus on the scaling limits. We give formal conjectures
and discuss them based on the analysis of the data. The data will
be presented in the form of figures. In the sequel, the jump rate
$\gamma$ of the RE plays a central role. Throughout the paper we
assumed, for simplicity, the RW jumping at rate $1$, if instead we
let it jump at any other rate $\lambda>0$, the same results would
hold replacing $\gamma$ by $\gamma/\lambda$.

\subsection{Asymptotic speed}
\label{S3.1}

Denote by $X_n=X_n(\rho,p,\gamma)$ the position of the RW in
\eqref{rwdef} after $n$ exponential times of rate $1$. In this
section we analyze the behavior of the asymptotic speed which we
obtained by evaluating for large $n$

\be{SimSpeed}v_n=v_n(\rho,p,\gamma):=\frac{1}{nM}\sum_{i=1}^{M}X_n^{(i)},\ee
over a sample of $M$ independent experiments (the values of $M$'s and $n$'s
will be specified in the figures.)

\bcj{Existence} Let $\gamma>0$, $(p,\rho)\in[1/2,1)\times[1/2,1)$,
and assume $\xi$ is the SSE, then $P^\mu-\xi$ a.s.
$$\exists \lim_{t\rightarrow\infty}\frac{X_t(p,\rho,\gamma)}{t}=:v(p,\rho,\gamma)\in\R.
$$\ecj

Conjecture \ref{Existence} should hold in great generality at
least for translation invariant RE. At the present state of the
art, the existence of an almost sure constant speed has been
proven for dynamic REs with ``good'' mixing properties in space
and time (see \cite{AvdHoRe11,BaZe06,dHoSaSi11,ReVo11}) except in
\cite{AvSaVo11} which instead uses a strong elliptic condition. In
particular, Conjecture \ref{Existence} is a rigorous statement if
$\xi$ is the ISF (see \cite{AvdHoRe11,ReVo11}).

\bcj{rhoV}Let $\xi$ be either the SSE or the ISF. For any
$\gamma>0$, $p\in(1/2,1)$, the function $\rho \longmapsto
v(p,\rho,\gamma)$ is continuous and non-decreasing. \ecj

Note that the monotonicity is trivial once the existence of
$v(p,\rho,\gamma)$ is given. Indeed, it follows by the fact that
for any $\rho<\rho'$, the RWs $X_t(p,\rho,\gamma)$ and
$X_t(p,\rho',\gamma)$ can be coupled so that they remain ordered.
Figure \ref{speedp} below illustrates the mentioned monotonicity.
In particular, it refers to the outcome of the numerics produced
in the case of the SSE with $p=0.8$. We remark that the same
qualitative picture holds in the ISF case and for any other choice
of $p\in(1/2,1)$.

\newpage
\begin{figure}[hbtp]
\begin{center}
\includegraphics[width=.8\textwidth]{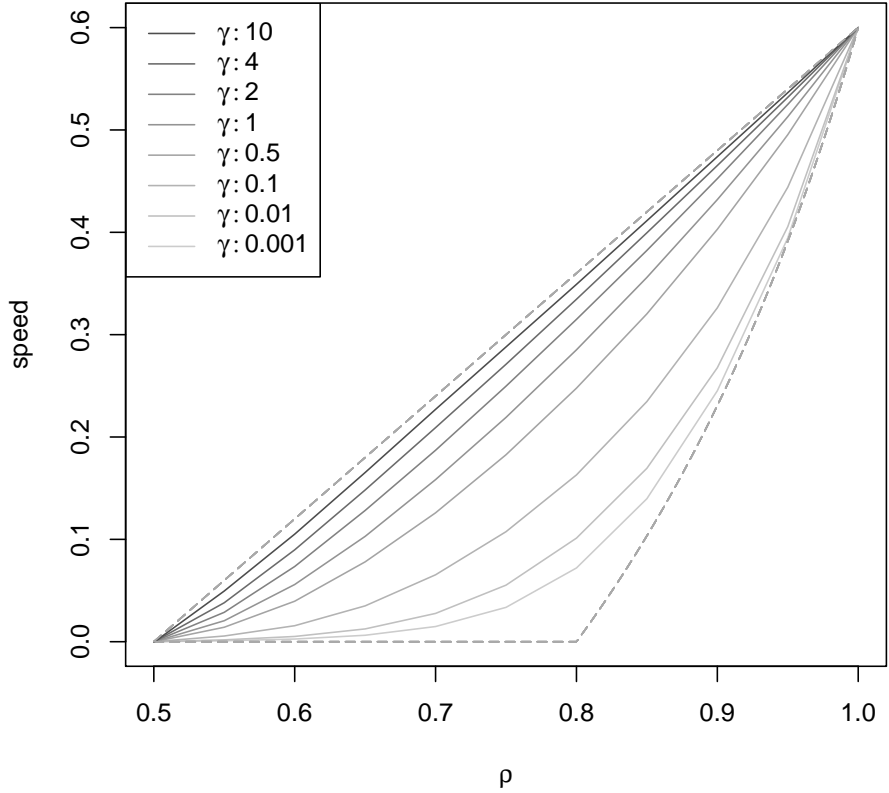}
\end{center}
\caption{The lower and the upper dashed curves correspond to the
speeds in the static \eqref{StSpeed1} and the averaged medium
\eqref{AvMedSpeed} cases, respectively. The different solid curves
represent the function $\rho \longmapsto v(0.8,\rho,\gamma)$ at
the different $\gamma$'s specified. In particular, each curve is
obtained after interpolating $11$ points at distance $0.05$, where
each point has been obtained by averaging over a sample with at least
$M=O(10^3)$ simulations of the empirical speed $v_n$ in
\eqref{SimSpeed} with $n=2^{16}$. This picture has been produced
by using the SSE as RE.} \label{speedp}
\end{figure}

\bcj{gammaV}Let $\xi$ be either the SSE or the ISF. Then the
function $\gamma \longmapsto v(p,\rho,\gamma)$ is continuous and
non-decreasing. Moreover \be{0Continuity}\lim_{\gamma\downarrow
0}v(p,\rho,\gamma)=v_{\mathrm{static}}(p,\rho),\ee
\be{InfinityContinuity}\lim_{\gamma\rightarrow
\infty}v(p,\rho,\gamma)=v_{\mathrm{averaged}}(p,\rho).\ee \ecj

\begin{figure}[hbtp]
\begin{center}
\includegraphics[width=.8\textwidth]{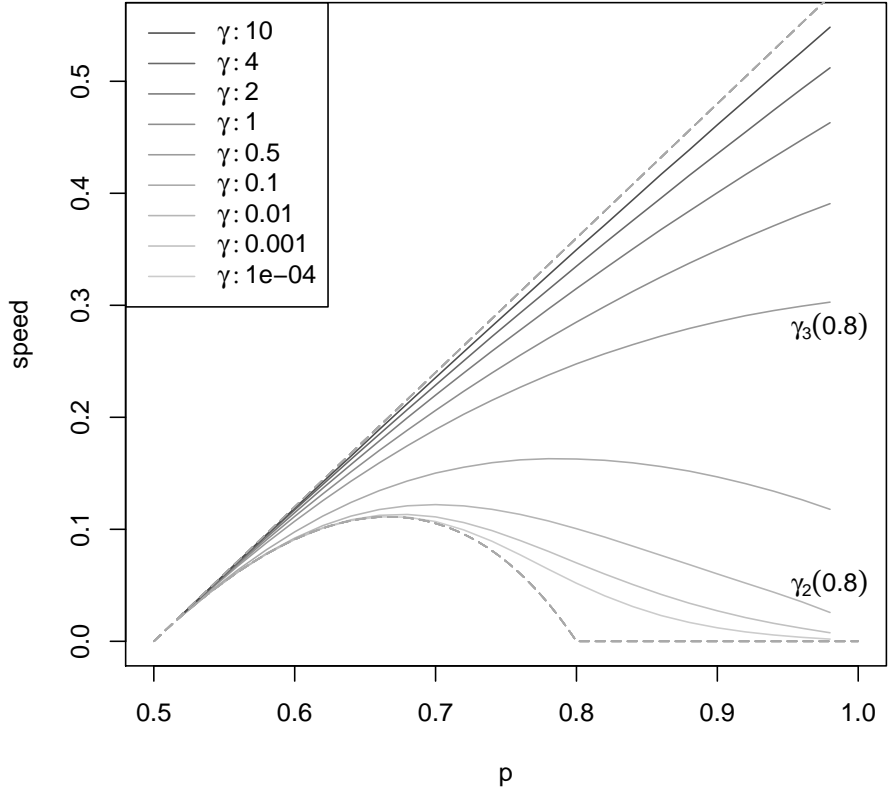}
\end{center}
\caption{$\rho=0.8$. The lower and the upper dashed curves
correspond to the speeds in the static \eqref{StSpeed1} and the
averaged medium \eqref{AvMedSpeed} cases, respectively. The
different solid curves represent the function $p \longmapsto
v(p,0.8,\gamma)$ at the different specified $\gamma$'s. Each solid
curve is obtained after interpolating $26$ points at distance
$0.02$, where each point has been obtained by averaging over the
outcome of samples with at least $M=O(10^3)$ independent simulations of the
empirical speed $v_n$ in \eqref{SimSpeed} with $n\geq 2^{18}$.
This plot has been produced by using the SSE as a RE.}
\label{speedrho}
\end{figure}

Conjecture \ref{gammaV} is supported by the results of the
experiments presented in Figures \ref{speedp} and \ref{speedrho}.
To state the next conjecture which describes the behavior of the
speed as a function of $p$, we introduce some quantities defined
in terms of $\gamma$ for fixed $\rho\in(1/2,1)$.

\bcj{pV} Let $\xi$ be either the SSE or the ISF. Fix
$\rho\in(1/2,1)$. Define \be{gamma1}\gamma_1(\rho):=\inf\{\gamma>0
:v(p,\rho,\gamma)>0 \text{ for all } p> 1/2\},\ee

\be{gamma2}\begin{aligned}\gamma_2(\rho)&:=\inf\{\gamma>0 :
p \longmapsto v(p,\rho,\gamma) \text{ is concave } \},\end{aligned}\ee

\be{gamma3}\gamma_3(\rho):=\inf\{\gamma>0 : p \longmapsto
v(p,\rho,\gamma) \text{ is non-decreasing } \}.\ee

Then, the function $p \longmapsto v(p,\rho,\gamma)$ is continuous.
Moreover, \be{gammas}0\leq
\gamma_1(\rho)<\gamma_2(\rho)<\gamma_3(\rho)<\infty.\ee\ecj

Conjecture \ref{pV} states the existence of several critical
$\gamma$'s for which we see a different behavior of the speed as a
function of $p$. Figure \ref{speedrho} below shows such a scenario
in the SSE case at a given $\rho=0.8$. Again, the same qualitative
picture holds for any other $\rho\in(1/2,1)$ or by considering the
ISF. Above $\gamma_3(\rho)$ the function $p \longmapsto
v(p,\rho,\gamma)$ is increasing, while it starts to become
non-monotone for $\gamma<\gamma_3(\rho)$. Below $\gamma_2(\rho)$
it looses the concavity and for $\gamma\leq\gamma_1(\rho)$ it
possibly starts to have a vanishing piece. A crucial issue is to
understand if

\be{OpenPb}\gamma_1(\rho)>0.\ee

The positivity of $\gamma_1(\rho)$ would imply a surprising
transient regime with zero speed which so far has been proven only
in the static case (see Section \ref{S2.1}). In the case of the
ISF, it follows from \cite{ReVo11} that $\gamma_1(\rho)=0$. It
might still be that $\gamma_1(\rho)>0$ in the SSE case.
Unfortunately, we feel not able to conjecture anything based on
our numerics since on one hand in the corresponding region
$\gamma\ll 1$, at the time scales we could achieve, the simulation
may just be a weak perturbation of the static case, and on the
other hand this is a very delicate phenomenon to test
statistically since of course in any transient regime the expected
speed at finite time is strictly positive. We therefore leave
\eqref{OpenPb} as a key open problem of this model.

The loss of monotonicity for low $\gamma$'s is related to the
strength of the traps. In the static case, the speed $p
\longmapsto v(p,\rho,\gamma)$ looks like the dashed lower curve in
figure \ref{speedrho} at any fixed $\rho\in(1/2,1)$. For $p$ big
enough it starts to become decreasing until it vanishes.
Intuitively, this is saying that when we increase $p$, no matter
what the size of a typical trap is, the holes tend to act almost
as reflecting barriers.

\begin{figure}[hbtp]
\begin{center}
\includegraphics[width=.46\textwidth]{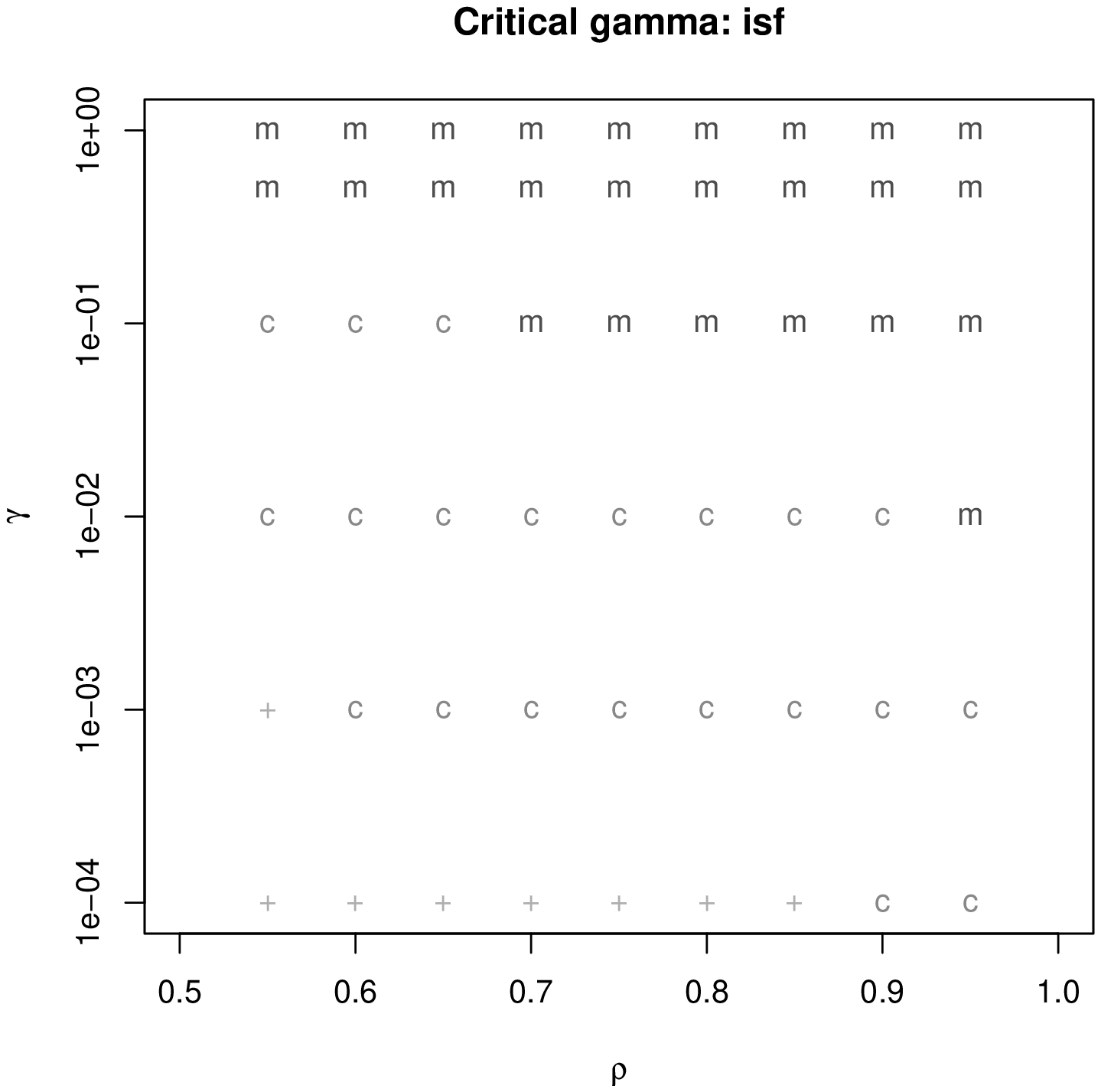}\\
\includegraphics[width=.46\textwidth]{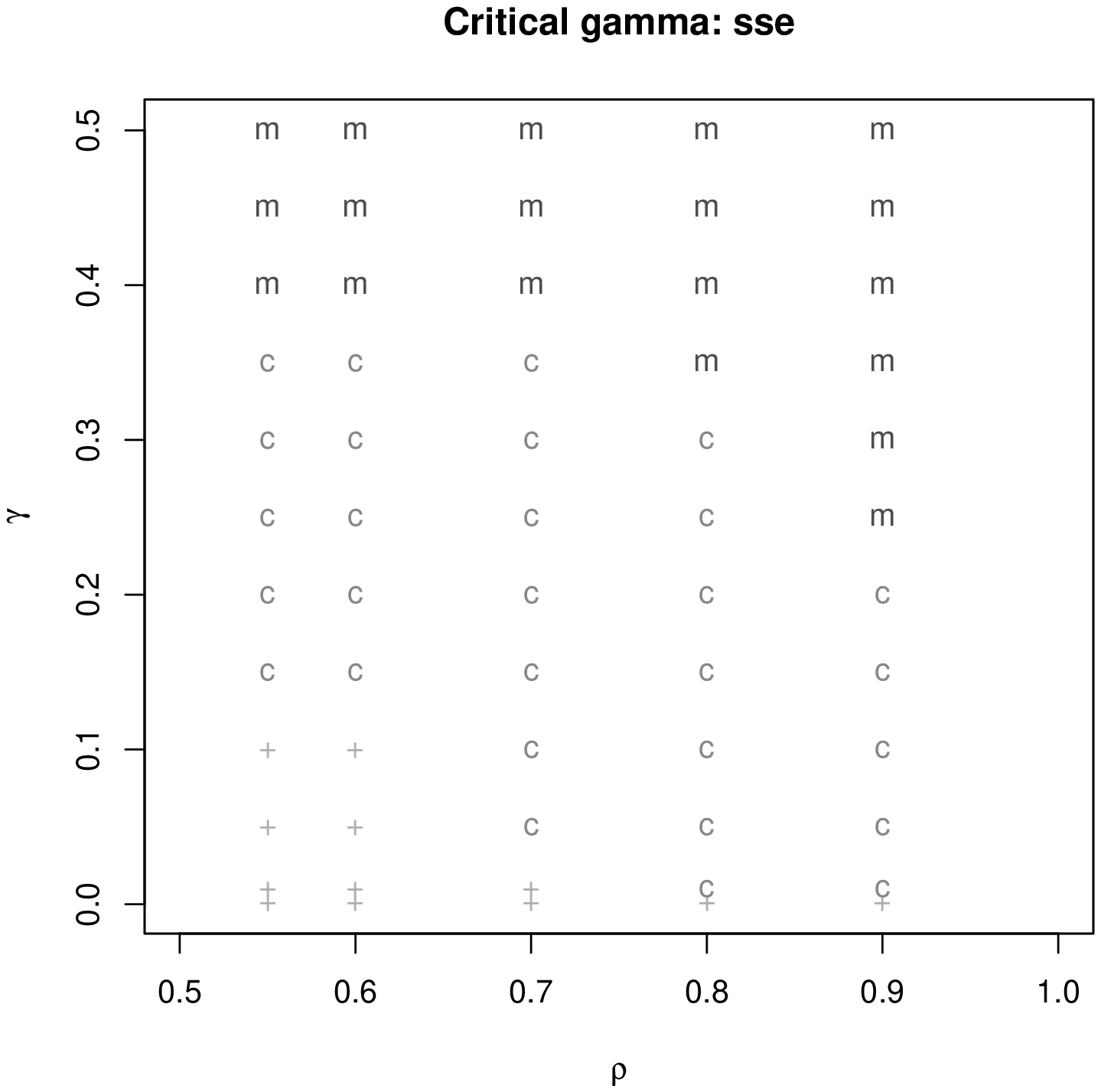}
\end{center}
\caption{Quantitative phase diagrams describing Conjecture
\ref{pV} in the ISF and the SSE environments. To each
$(\rho,\gamma)$ we associate the symbols m, c or +. They mean that
the corresponding speed is a monotone, concave and non-concave
function of $p$. $\gamma_3(\rho)$ would correspond to the curve in
between the regions of m's and c's, while $\gamma_2(\rho)$ would
be the line separating the regions with c's and +'s.}
\label{criticalgammas}
\end{figure}

\newpage
Figure \ref{criticalgammas} presents a quantitative version of
Conjecture \ref{pV} in both the ISF and SSE cases.

\newpage
\subsection{Scaling limits and phase diagram: SSE}
\label{S3.2} We now turn to the analysis of the scaling limit of
$X$. The ISF case belongs to the class of dynamic RE with strong
mixing properties and it is known that for such dynamics for any
$\gamma>0$, $X$ satisfies a functional CLT (see e.g.
\cite{ReVo11}). In the sequel we therefore focus only on the SSE
case. We just mention that we checked numerically the known
diffusivity of the ISF case. We obtained excellent agreement in
this case except for very few points close to degenerate cases,
i.e.\ for $\gamma$ too small and $(p,\rho)\approx(1,1)$ (see also
the discussion in Section \ref{S4.2}).

Our main conjecture is:

\bcj{Scaling} Let $\xi$ be the SSE. There exist two monotone
functions $\widetilde{\gamma_1},\widetilde{\gamma_2}:
[1/2,1)^2\longmapsto\R^+$ with

\begin{enumerate}
\item[i)] $\widetilde{\gamma_i}$ equals zero only on the sets
$A_i, i=1,2$, defined by \be{A1}A_1:=\{(p,\rho)\in[1/2,1]^2:
\rho\leq f_1(p)\},\ee \be{A2}A_2:=\{(p,\rho)\in[1/2,1]^2: \rho\leq
f_2(p)\},\ee with $f_1(p)=\frac{(1-p)^2}{p^2}$ and
$f_2(p)=\frac{p-\sqrt{p(1-p)}}{2p-1}$(see Figure
\ref{ScalingStatic}), \item[ii)]
$\widetilde{\gamma_1}<\widetilde{\gamma_2}$ whenever they are
non-zero, \item[iii)] $\widetilde{\gamma_i}$ is increasing in $p$
and decreasing in $\rho$, for any $(p,\rho)\in[1/2,1)^2\setminus
A_i, i=1,2,$
\end{enumerate} and such that for any
$(p,\rho)\in[1/2,1)^2$, we have the following cases

\begin{enumerate}
\item$X(p,\rho,\gamma) \text{ is sub-diffusive for any }
\gamma<\widetilde{\gamma_1}(p,\rho)$ \item$X(p,\rho,\gamma) \text{
is super-diffusive for any }
\widetilde{\gamma_1}(p,\rho)<\gamma<\widetilde{\gamma_2}(p,\rho)$
\item$X(p,\rho,\gamma) \text{ is diffusive for any }
\gamma>\widetilde{\gamma_2}(p,\rho)$ or at
$\gamma=\widetilde{\gamma_1}(p,\rho)$
\end{enumerate} \ecj

\begin{figure}[hbtp]
\vspace{0.5cm}
\begin{center}
\includegraphics[width=.45\textwidth]{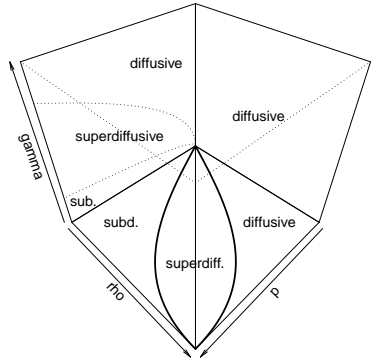}
\end{center}
\caption{\small A qualitative picture of the phase diagram
described in Conjecture \ref{Scaling} in the degenerate cases
$\rho=1/2$ (recurrent case) and $p=1/2$ (SSRW). The lower and
upper dotted curves at $\rho=1/2$ represent
$\widetilde{\gamma_1}(p,1/2)$ and $\widetilde{\gamma_2}(p,1/2)$,
respectively. } \label{cuberec}
\end{figure}

\begin{figure}[hbtp]
\vspace{0.5cm}
\begin{center}
\end{center}
\begin{center}
\includegraphics[width=.8\textwidth]{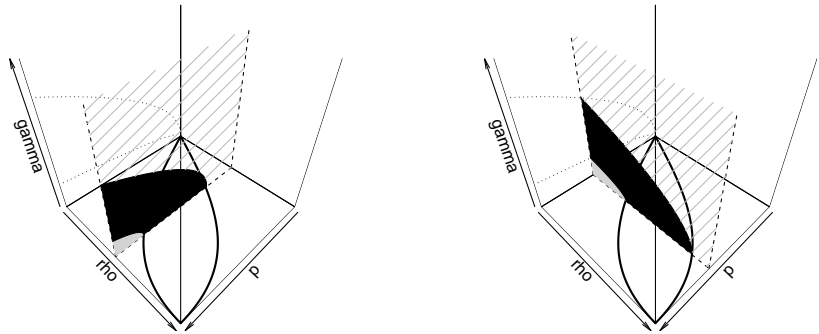}
\end{center}
\caption{\small Two qualitative sections of the phase diagram
described in Conjecture \ref{Scaling} for fixed $\rho$ (on the
left) and $p$ (on the right). The areas of those sections in gray,
black and striped correspond to sub-diffusive, super-diffusive and
diffusive regimes, respectively. In particular, the lower and
upper curves bounding the black areas represent
$\widetilde{\gamma_1}$ and $\widetilde{\gamma_2}$, respectively.}
\label{cubesect}
\end{figure}

This conjecture is the most interesting novel result of our
numerics. Figures \ref{cuberec} and \ref{cubesect} show the
qualitative scenario stated in Conjecture \ref{Scaling}. Recall
the ``super-diffusive leaf'' in Figure \ref{ScalingStatic}, note
that the functions $f_1$ and $f_2$ represent the lower and the
upper boundary of the leaf, respectively. The fact that for the
ISF (and more generally for RE with space-time correlation with
exponential or fast polynomial decay) at any $\gamma>0$ we have a
diffusive scaling, can be rephrased by saying that for any
$\gamma>0$ the leaf vanishes and any point $(p,\rho)\in[1/2,1)^2$
corresponds to diffusive regime. On the other hand, in case of the
SSE, Conjecture \ref{Scaling} says that as soon as we switch on
the SSE dynamics (i.e. for small $\gamma>0$) the leaf is still
present and starts to move towards the $p$-axis as $\gamma$
increases until a certain critical $\gamma$ for which the leaf
completely disappears. Note that the disappearance of the leaf for
$\gamma$ big enough is consistent with the fact that as $\gamma$
increases we get closer and closer to the averaged medium case
which is clearly diffusive. This observed phenomenon (although it
could still be local, see Section \ref{S4.2}) suggests that due to
the slow-mixing properties of the exclusion dynamics, traps play a
crucial role to determine the scaling limit of $X$. In particular,
depending on the ratio of the jump rate of the walker and the one
of the SSE, we can observe diffusivity or not.

\begin{figure}[hbtp]
\vspace{0.5cm}
\begin{center}
\includegraphics[width=\textwidth]{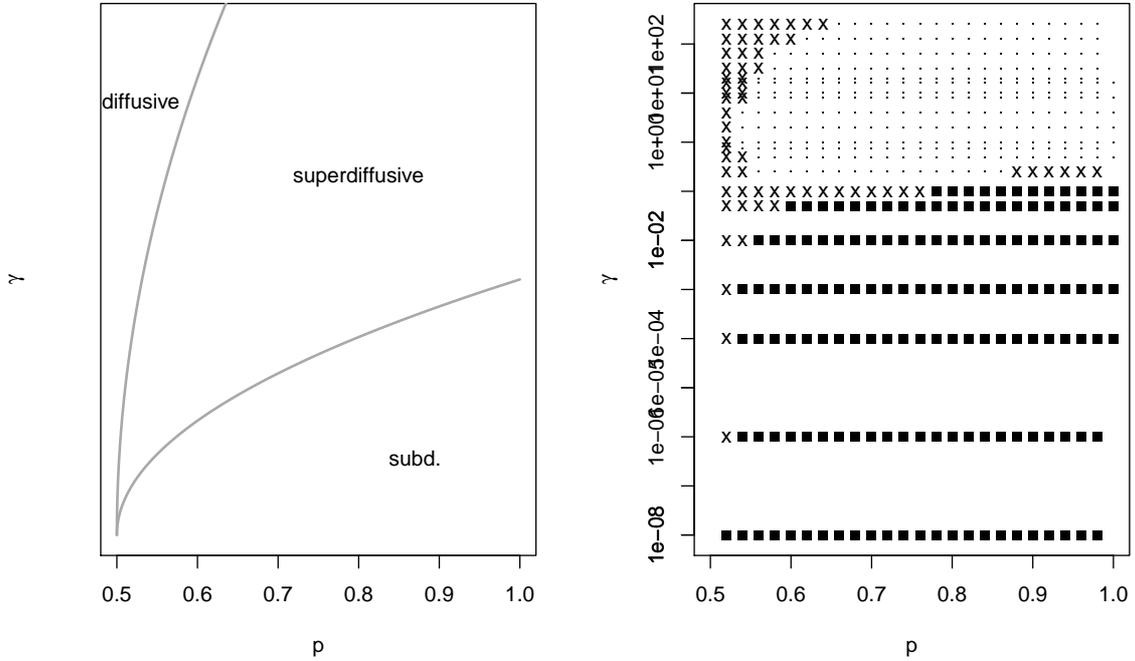}
\end{center}
\caption{\small The section of the phase diagram in the recurrent
case, i.e. $\rho=0.5$. On the left the qualitative picture, on the
right the outcome of our experiments supporting the qualitative
phase diagram. The crosses, the dots, and the black squares mean
that for the corresponding $(\rho,\gamma)$ points our test gave an
exponent equal (diffusive), bigger (super-diffusive), or smaller
(sub-diffusive) than $1/2$, respectively. More precisely, we
assigned a cross at any estimated exponent in between
$[0.49,0.51]$. Note that as in the lower part of the leaf for the
static case (see Figure \ref{ScalingStatic}), on the line dividing
the sub- and super-diffusive regimes the scaling is
diffusive.}\label{rainbowrec}
\end{figure}

\begin{figure}[hbtp]
\vspace{0.5cm}
\begin{center}
\includegraphics[width=.7\textwidth]{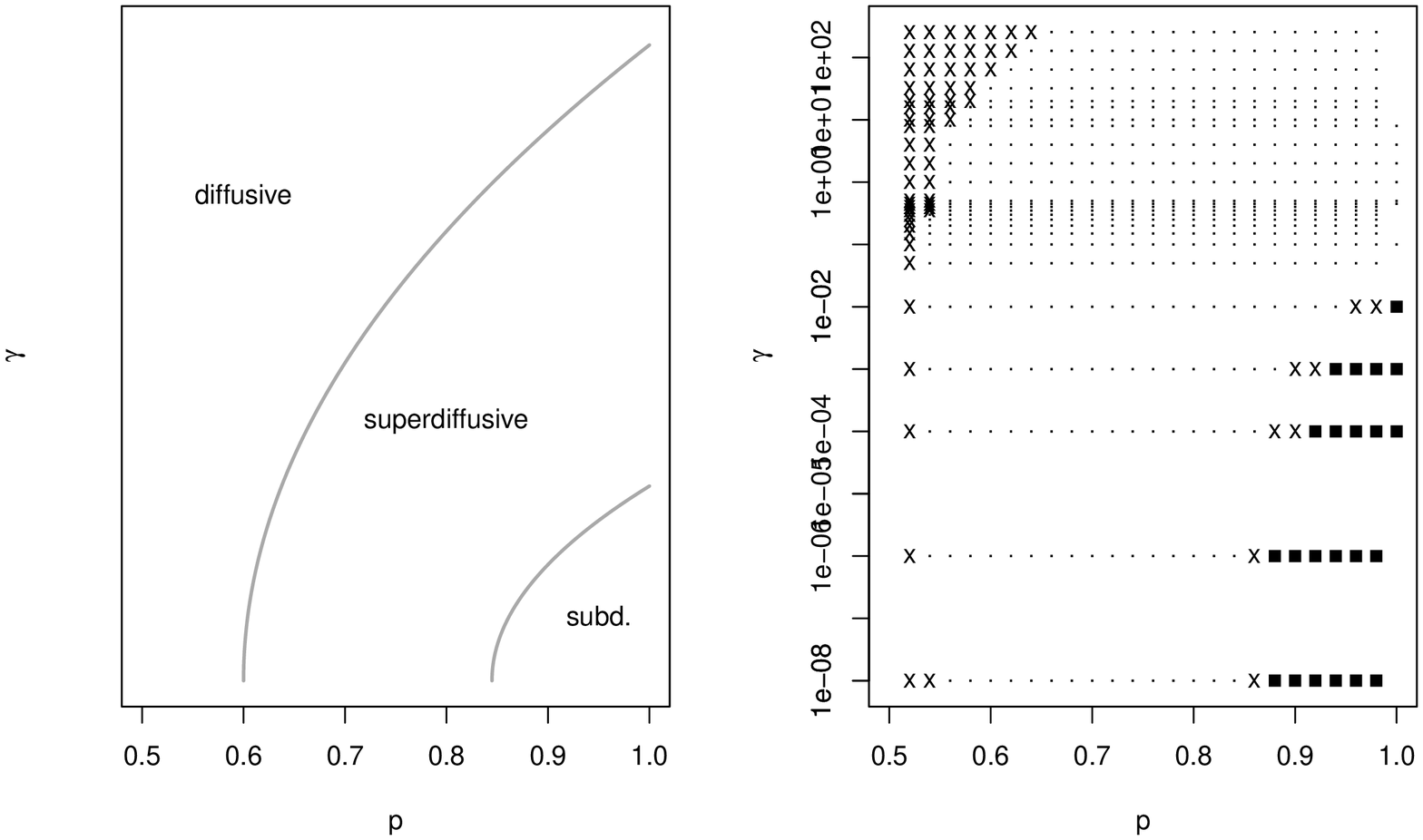}
\end{center}
\caption{\small The section of the phase diagram at $\rho=0.7$. On
the left the qualitative picture, On the right the outcome of our
experiments. As in the previous picture, crosses, dots, and
squares correspond to diffusive, super- and sub-diffusive
regimes, respectively.} \label{rainbowrho}
\end{figure}

\begin{figure}[hbtp]
\vspace{0.5cm}
\begin{center}
\includegraphics[width=.7\textwidth]{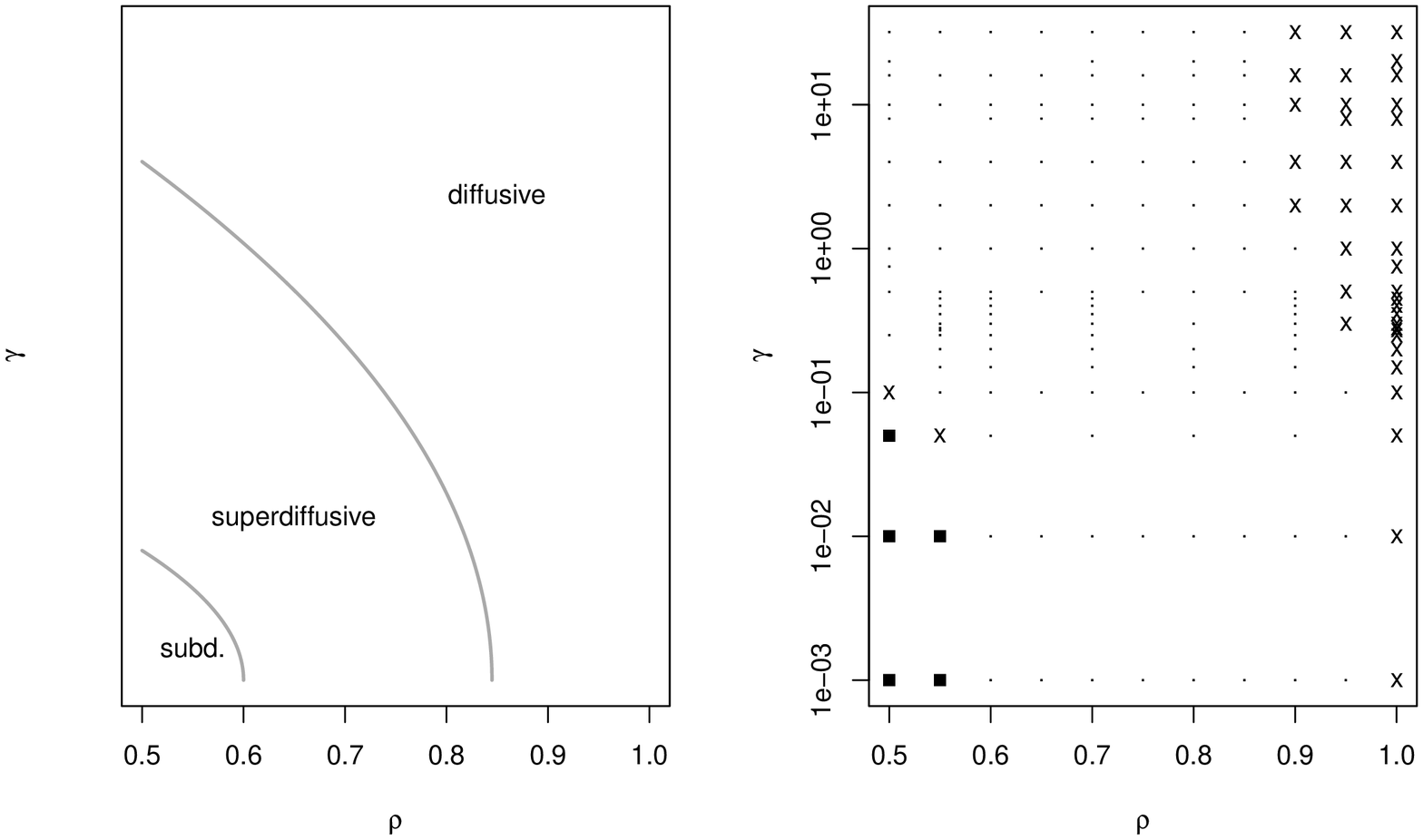}
\end{center}
\caption{\small The section of the phase diagram at $p=0.7$. On
the left the qualitative picture, on the right the outcome of our
experiments. As in the previous pictures, crosses, dots, and
squares correspond to diffusive, super- and sub-diffusive
regimes, respectively.} \label{rainbowp}
\end{figure}

\begin{figure}[hbtp]
\vspace{0.5cm}
\begin{center}
\includegraphics[width=\textwidth]{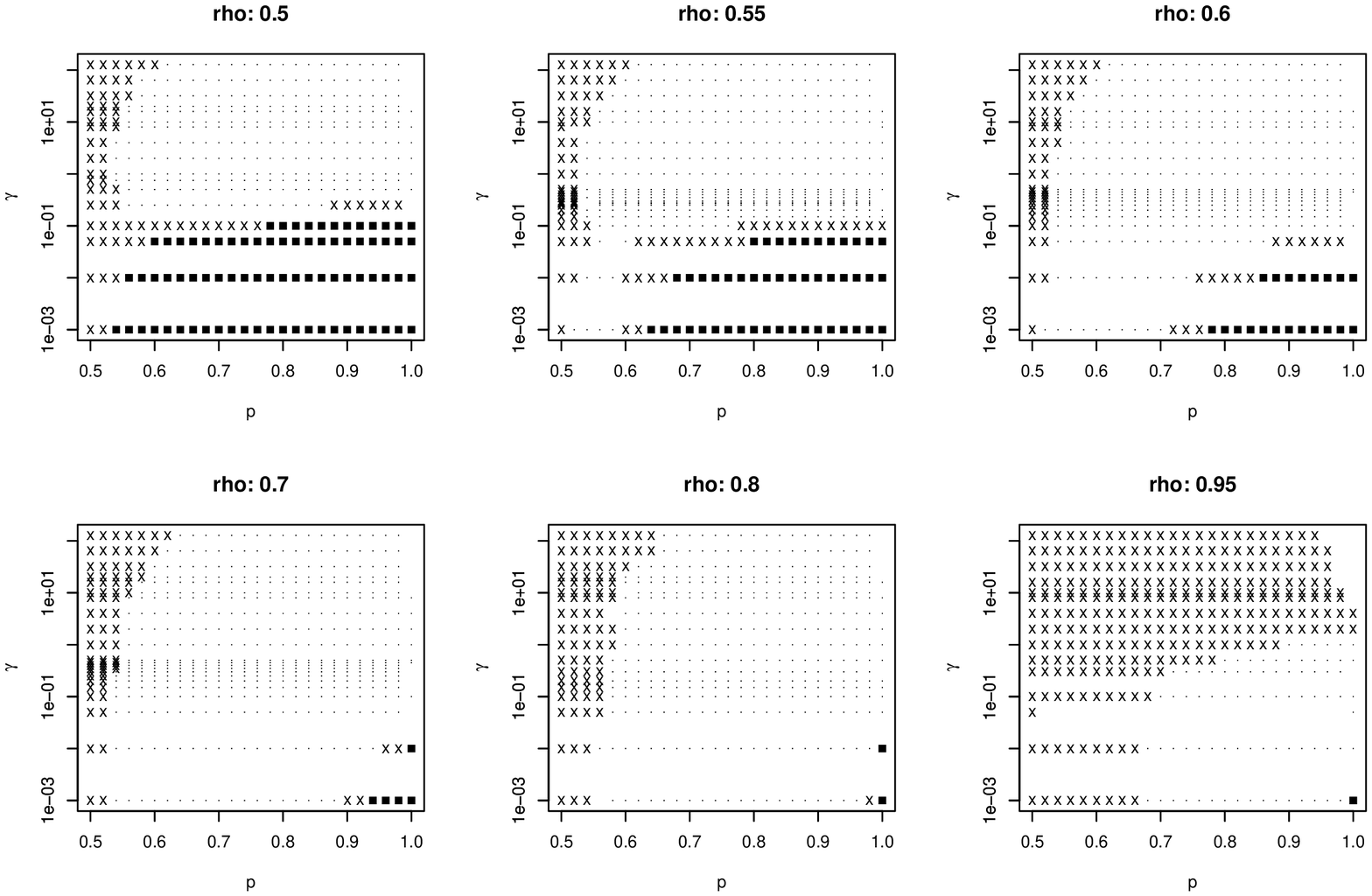}
\end{center}
\caption{\small Sections of data supporting the phase diagram in
Conjecture \ref{Scaling} increasing in $\rho$. Crosses, dots, and
squares correspond to diffusive, super- and sub-diffusive
regimes, respectively.}\label{rhosequence}
\end{figure}

\begin{figure}[hbtp]
\vspace{0.5cm}
\begin{center}
\includegraphics[width=\textwidth]{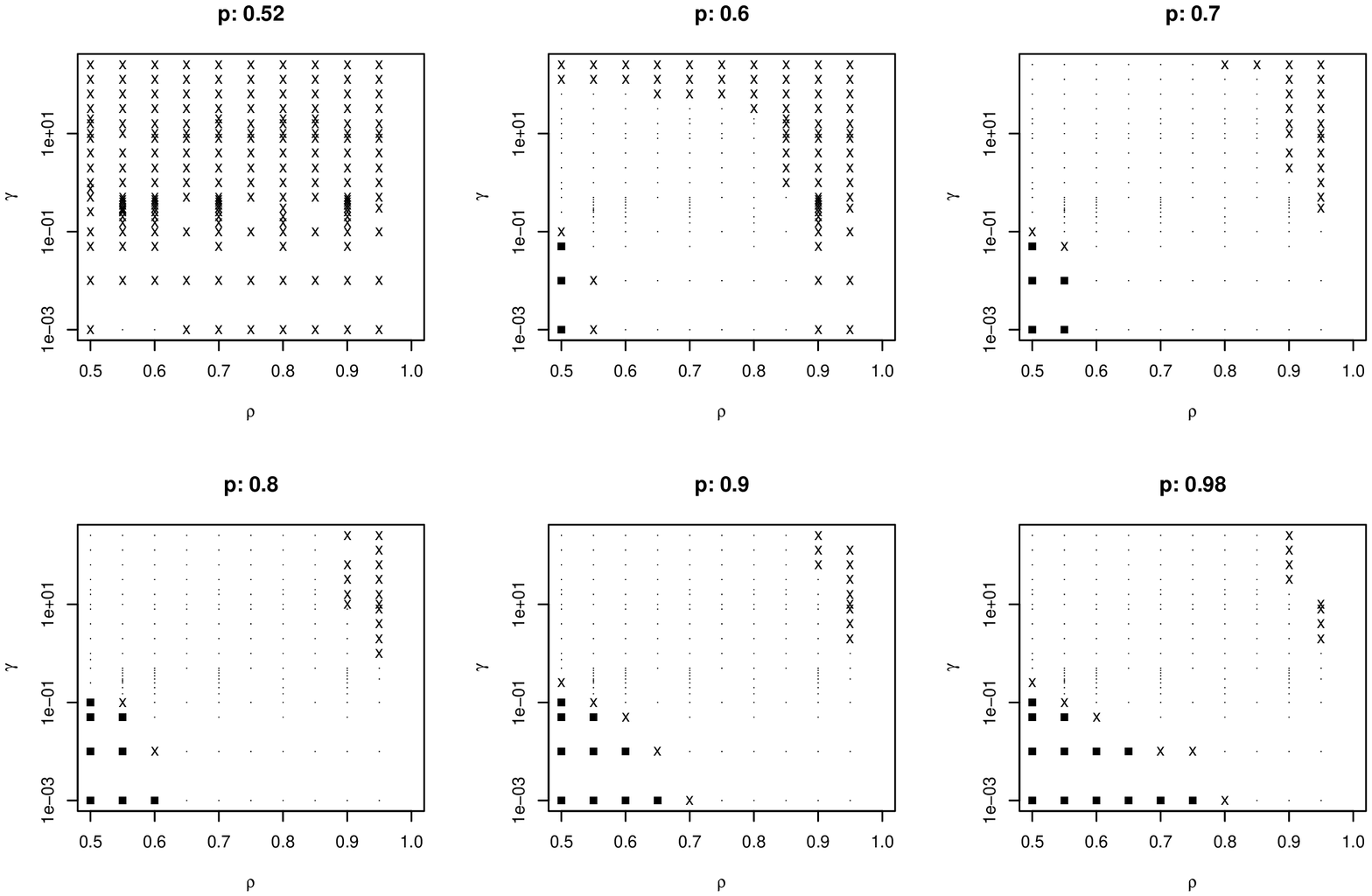}
\end{center}
\caption{\small Sections of data supporting the phase diagram in
Conjecture \ref{Scaling}  increasing in $p$. Crosses, dots, and
squares correspond to diffusive, super- and sub-diffusive
regimes, respectively.}\label{psequence}
\end{figure}

\newpage
Figures
\ref{rainbowrec}--\ref{rainbowrho}--\ref{rainbowp}--\ref{rhosequence}--\ref{psequence}
show some of the data supporting Conjecture \ref{Scaling}. In
Section \ref{S3.3} we describe how we obtained these phase
diagrams.

\newpage
\subsection{Scaling exponents: estimating the variance for SSE}
\label{S3.3}

This section is devoted to the description of the estimates we used to
obtain the phase diagram described in Conjecture \ref{Scaling},
i.e.\ the scaling exponents of $X$ in the SSE environment.

Figure \ref{explicitsect} below shows the same section as in
Figure \ref{rainbowrec} where instead of marking ``diffusive and
non-diffusive points'', we give the explicit values of the scaling
exponents.

\begin{figure}[hbtp]
\vspace{0.5cm}
\begin{center}
\includegraphics[width=\textwidth]{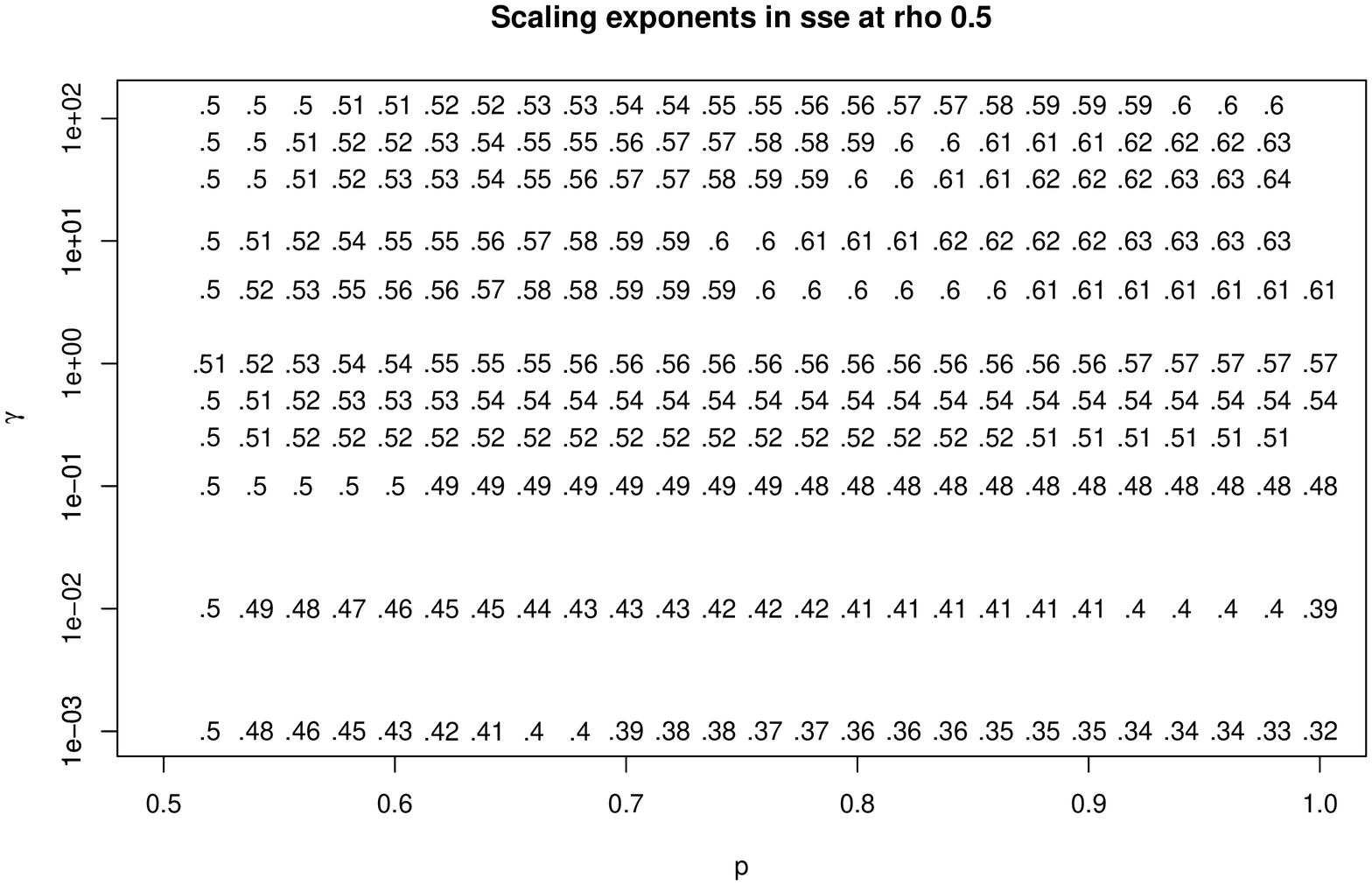}
\end{center}
\caption{\small A quantitative version of the phase diagram in
Figure \ref{rainbowrec}. The numbers written at each point
$(p,0.5,\gamma)$ are the scaling exponents
$\alpha^\star(p,\rho,\gamma)$ that were estimated as described
below.} \label{explicitsect}
\end{figure}

\begin{figure}[hbtp]
\begin{center}
\includegraphics[width=\textwidth]{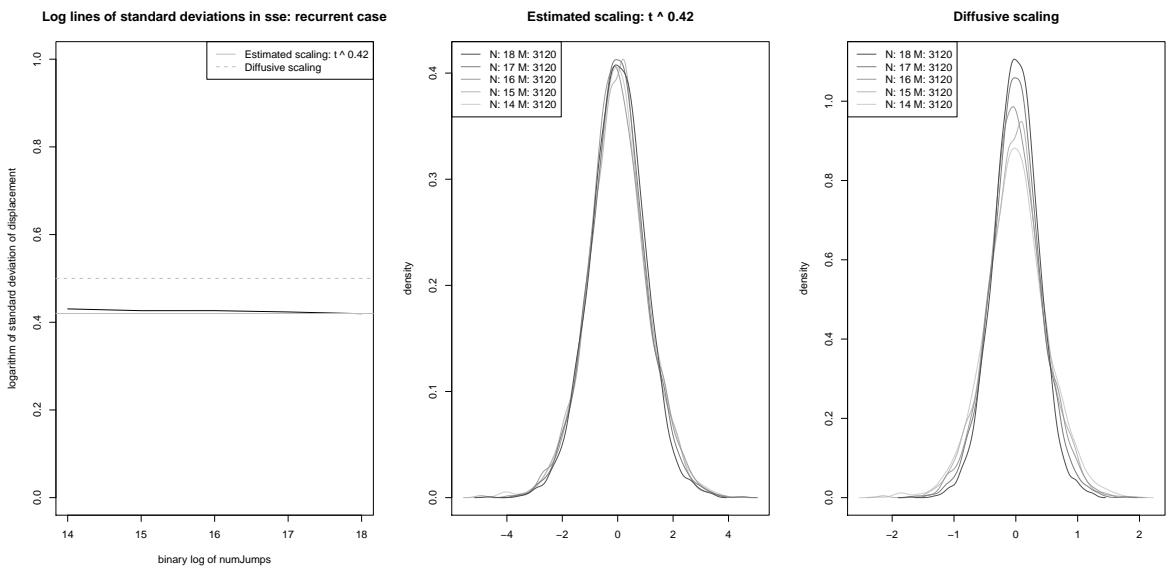}
\end{center}
\caption{\small Recurrent, sub-diffusive.
$(p,\rho,\gamma)=(0.76,0.5,0.01)$. The left picture illustrates
the behavior of the logarithm of the standard deviation of $X_n$
(the solid line). The dotted line represents the estimated
exponent $\alpha^\star(0.76,0.5,0.01)$, while the dashed one
represents the diffusive value $0.5$. The middle plot shows that
the densities obtained at different times $n=2^N$ rescaled by
$n^{\alpha^\star}$ almost perfectly match. In the right most plot,
the same densities under diffusive scaling, looking carefully at
the tails and around $0$, one can see that they are ordered in
time and tend to concentrate in $0$ as time increases. $M$ denotes
the sample size. The data at different times $n$ were obtained
with independent experiments.} \label{rsub}
\end{figure}

\begin{figure}[hbtp]
\vspace{0.5cm}
\begin{center}
\includegraphics[width=.65\textwidth]{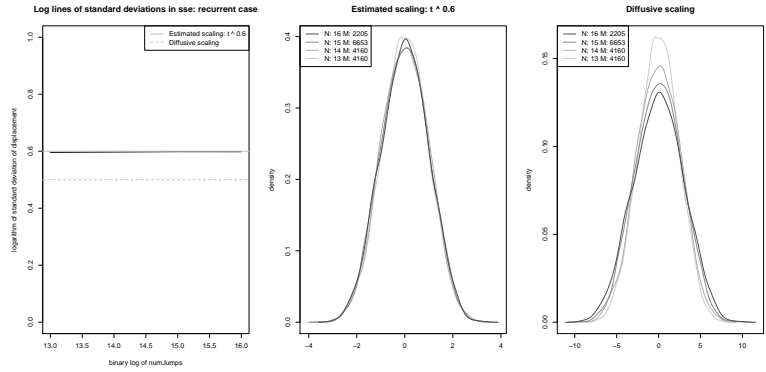}
\end{center}
\caption{\small Recurrent, super-diffusive.
$(p,\rho,\gamma)=(0.8,0.5,4)$. As in the previous figure, the left
picture illustrates the behavior of the logarithm of the standard
deviation of $X_n$ (the solid line). The dotted line represents
the estimated exponent $\alpha^\star(0.8,0.5,4)$, while the dashed
one represents the diffusive value $0.5$. The middle plot shows
that the densities obtained at different times $n=2^N$ rescaled by
$n^{\alpha^\star}$ almost perfectly match. In the right most plot,
the same densities under diffusive scaling, note that in this case
they are ordered in time and tend to vanish as time increases. $M$
denotes the sample size. The data at different times $n$ were
obtained with independent experiments.} \label{rsup}
\end{figure}

\begin{figure}[hbtp]
\begin{center}
\includegraphics[width=.65\textwidth]{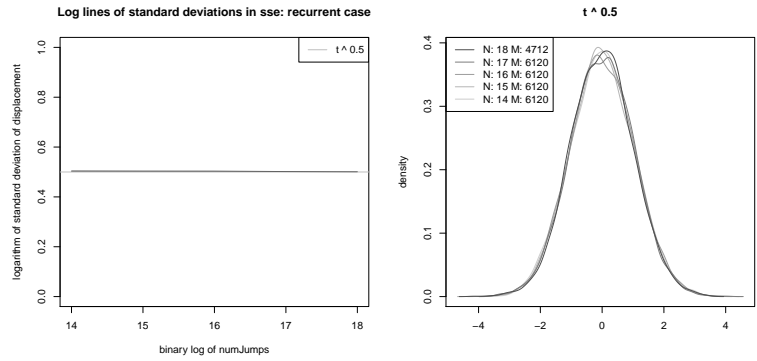}
\end{center}
\caption{\small Recurrent, diffusive.
$(p,\rho,\gamma)=(0.56,0.5,0.1)$. As in the previous figures, the
left picture illustrates the behavior of the logarithm of the
standard deviation of $X_n$. The dotted line represents the
estimated exponent $\alpha^\star(0.6,0.5,0.1)$. The right plot
shows that the densities obtained at different times $n=2^N$
coincide under diffusive rescaling.} \label{rdif}
\end{figure}

\begin{figure}[hbtp]
\begin{center}
\includegraphics[width=.8\textwidth]{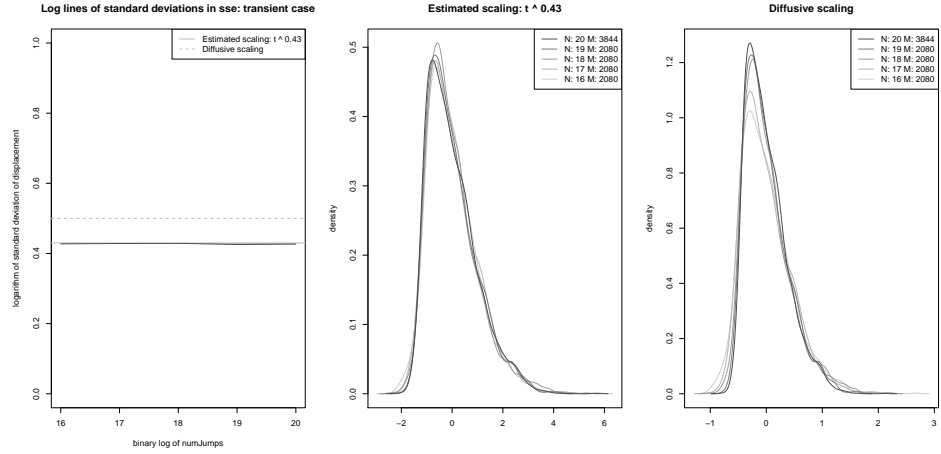}
\end{center}
\caption{\small Transient, sub-diffusive. $(p,\rho,\gamma)=(0.74,0.55,0.001)$.
Estimated reliable time: $\bar{n}=2^{38}$ (see Section \ref{S4.2}).}
\label{tsub}
\end{figure}

\begin{figure}[hbtp]
\begin{center}
\includegraphics[width=.8\textwidth]{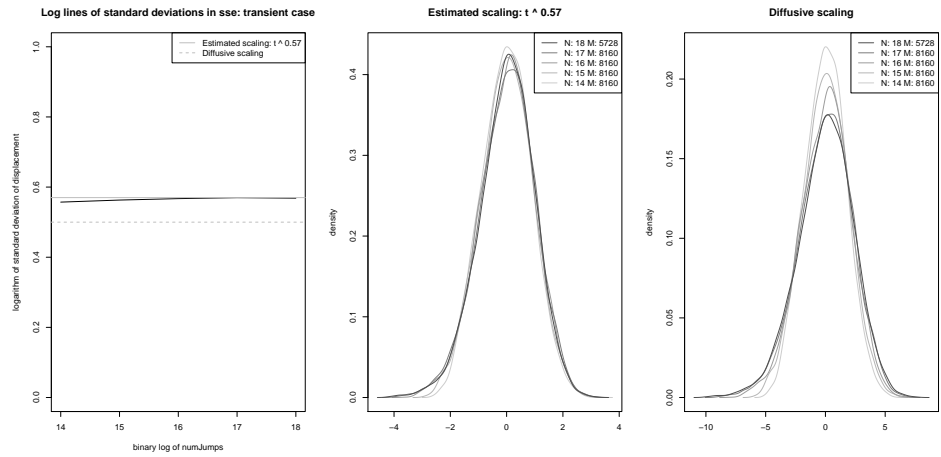}
\end{center}
\caption{\small Transient, super-diffusive.
$(p,\rho,\gamma)=(0.86,0.6,0.25)$. Estimated reliable time:
$\bar{n}=2^{12}$, effective running time $n=2^{18}$ (see Section
\ref{S4.2}).} \label{tsup}
\end{figure}

\begin{figure}[hbtp]
\begin{center}
\includegraphics[width=.6\textwidth]{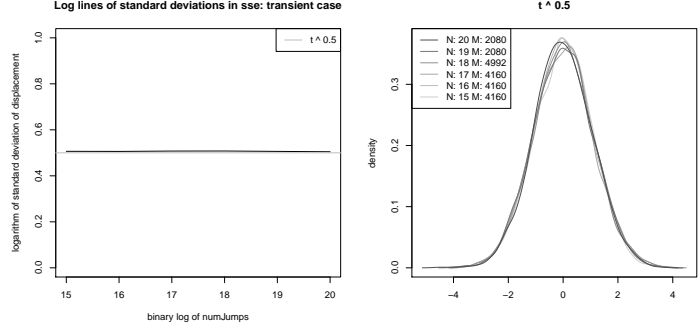}
\end{center}
\caption{\small Transient, diffusive.
$(p,\rho,\gamma)=(0.56,0.9,0.001)$. Estimated reliable time:
$\bar{n}=2^{305}$ (see Section \ref{S4.2}).} \label{tdif}
\end{figure}

Each of these exponents is obtained by analyzing the variances of
$X$ at different times. More precisely, for each fixed triple
$(p,\rho,\gamma)$, for fixed number of jumps $n$, we compute the
sample standard deviation of $X_n(p,\rho,\gamma)$, namely,
$$SD_n:=\sqrt{\frac{1}{M-1}\sum_{i=1}^{M}
\left(X^{(i)}_n-v_n n\right)^2},$$ over a sample of $M$
independent experiments, with $M$ of order at least $10^3$ (the
values of $M$ are specified in the figures) and $v_n$ as in
\eqref{SimSpeed}. We then evaluate the function
$$\alpha(n):=\frac{\log (SD_n)}{\log n}$$ on a $\log n$ scale, see
Figures
\ref{rsub}--\ref{rsup}--\ref{rdif}--\ref{tsub}--\ref{tsup}--\ref{tdif}.
In particular, for each value of $n$ (in the experiments $n$ grows
like $2^N$, with $N$ specified in the figures) we have performed
independent experiments. The proper scaling exponent is
approximated by the value of $\alpha(n)$ with the biggest $n$
which we denote by $\alpha^\star= \alpha^\star(p,\rho,\gamma)$
(those are the numbers in Figure \ref {explicitsect}). Indeed,
assume $SD_n=c n^\alpha$, for some positive $c=c(p,\rho,\gamma)$,
then $\alpha(n)=\frac{\log c}{\log n}+\alpha$ converges to
$\alpha$ as $n$ goes to infinity.

The next figures show examples of these estimates in all described
regimes. Therein we plotted and overlapped the empirical densities
of $X_n-v_n n$ (obtained with samples of size $M$) for different
values of $n$ to see that indeed they coincide when rescaled by
$n^{\alpha^\star}$. In the non-diffusive cases (i.e.
$\alpha^\star\neq1/2$), we also add the plot of the same empirical
densities rescaled by $n^{1/2}$ to see that under diffusive
scaling they do not coincide. In particular, they are ordered in
time and the variances have the tendency of either vanishing
(sub-diffusive) or concentrating (super-diffusive).

\newpage

\newpage

\newpage

For the sake of completeness, we mention that we tried also some
other methods to guess the right scaling exponent. A quite apt but
not entirely stable method was to look for the exponent minimizing
the total variation distance between the densities at different
times. Another method was to fit an appropriate curve through the
standard deviations directly. All of them gave the same
qualitative picture we described here.

\section{Robustness of simulations and more general RWRE}
\label{S4}

\subsection{Algorithms}
\label{S4.1} To simulate the RW in the static RE, in the ISF and
in the SSE, we implemented an algorithm which can be briefly
described as follows.

\begin{enumerate}
\item[1)] Take as INPUT $(p,\rho,\gamma,
n)\in(0.5,1)\times[0.5,1)\times\R^+\times\N$, with $n$ being the
number of jumps of the RW $X$. \item[2)] Consider an interval $I$
of size $3n$ centered at the starting position of the RW.
\item[3)] Initialize the RE according to a Bernoulli product
measure of parameter $\rho$. \item[4)] At exponential time of rate
$1$, update the RW position according to the underlying state of
the RE. \item[5)] For each $x$ in $I$, at exponential time of rate
$\gamma$ update the RE at position $x$. \item[6)] Give as OUTPUT
the position of the RW after $n$ steps.
\end{enumerate}

Note that Step 5) depends on the RE. For the static case,
$\gamma=0$ implies that we never update the RE. In the ISF case,
due to the independence in space of the dynamics, we update the RE
only locally at the position of the RW $X$ (see Figure
\ref{Trajectories}). The case of the SSE has been implemented
using a version of the SSE either on the torus (i.e. with periodic
boundary conditions) or by sampling at rate $\gamma$ the state of
the RE at the boundaries of $I$ from $\nu_\rho$ (both approaches
produce the same outcome). In particular, the size of $I$ has been
chosen of size $3n$. Particles coming from the boundaries will
typically travel a distance of order $\sqrt{\gamma n}$ up to time
$n$ which guarantees that the dynamics of $X$ is not significantly
affected by those particles (note that the biggest $\gamma$ we
considered in our numerics is of order $100$).

\subsection{Heuristics versus locality}
\label{S4.2} Conjecture \ref{Scaling} suggests a completely
different scenario for the scaling limits depending on whether the
RE is \emph{fast} (ISF) or \emph{slowly mixing} (SSE). This latter
result is certainly the most interesting suggested by our
numerics. Indeed, in all the papers dealing with RWDRE only
diffusivity has been proven since, due to technical difficulties,
most of the tools available are not suitable to treat RE
presenting space-time correlations except under strong mixing
conditions. The non-diffusivity in the experiments for the SSE
could still be a local phenomenon. Namely, the time we run the
algorithm could in principle be too small and what we are
observing is only a perturbation of the static case or due to
other local effects: e.g. presence of constants, too small
$\gamma$, parameters in a neighborhood of the boundary of two
phases or close to degenerate cases. Therefore, to be more
confident that we are truly observing the right asymptotic
scalings, and to have some stronger evidence in favor of
Conjecture \ref{Scaling}, we need to understand whether the
running time of each experiment is large enough.

We present here a heuristic argument to guess a reliable running
time for small $\gamma$. The idea is that in each simulation the
RW should perform a number of jumps sufficiently large to ensure
that the time it takes the walker to cross a typical trap is
comparable with the time in which this trap dissolves.

Consider a trap of size $L$ of consecutive holes (a stretch of
size $L$ in $\Z$ where the state of the RE is $0$), assume the
walker starts at time $0$ on the left-most hole. Let $\tau_L$
denote the number of jumps the walker needs before reaching the
right-most hole of the trap. By using a standard gambler's ruin
argument we see that $$\E[\tau_L]\geq\frac{(p/q)^L-1}{(p/q)-1}.$$

On the other hand we can estimate the ``mixing time" of the trap
to be of order $(\gamma^{-1}L)^2$. Indeed, a particle at the
boundary of the trap (since it is performing a simple symmetric RW
at rate $\gamma$) would need in average this amount of time to
cross the trap.

If $\E[\tau_L]$ is much smaller than $(\gamma^{-1}L)^2$, this
would mean that within the time the RW crosses the trap, the
dissolvence effect due to the dynamics of the SSE is not big
enough to play a substantial role. Therefore we would like to
choose $L$ big enough so that the two quantities are at least
comparable, namely, we have to solve (when possible) the equation
$$\frac{(p/q)^L-1}{(p/q)-1}=(\gamma^{-1}L)^2.$$
Once we compute (numerically) the solution $L=L(p,\gamma)$ of this
equation, we have to be sure that the RW $X$ travels a distance
big enough so that it will meet with high probability such a big
trap. The probability of observing such a disaster can be easily
estimated by a geometric argument since the probability of a
stretch of hole of size $L$ is $(1-\rho)^L$. Therefore to be sure
that the RW $X(p,\rho,\gamma)$ will meet and cross at least one
disaster with high probability we have to run the algorithm for
$\bar{n}=\bar{n}(p,\rho,\gamma)$ steps, with $\bar{n}$ big enough
so that $X_{\bar{n}}\geq L(1-\rho)^L$ with high probability.

It turns out that the order of such a $\bar{n}(p,\rho,\gamma)$
varies a lot depending on $(p,\rho,\gamma)$. In most of the
parameter space is unfortunately too big to be achieved in a
reasonable amount of computing time. In the caption of Figures
\ref{tsub}-\ref{tsup}-\ref{tdif} we wrote the explicit
corresponding values of $\bar{n}$. The fact that in most of the
cases we could not run the algorithm up to such an $\bar{n}$, does
not imply that the outcome of the associated experiment is only
local and not asymptotically reliable, as the argument to deduce
$\bar{n}$ is only a rough estimate. On the other hand, the fact
that for several super-diffusive cases (as in the example in
Figure \ref{tsup}) we could achieve the associated $\bar{n}$ is a
strong suggestion in favor of Conjecture \ref{Scaling}.

\subsection{Concluding remarks}
\label{S4.3} As mentioned before, at the current state of the art,
most of the tools (regeneration, renormalization, coupling,
martingale approximation, etc.) developed to analyze RWRE (both
static and dynamic REs) are still inappropriate to deal with
space-time correlations unless the environment satisfies some
uniform and fast enough mixing condition. This is not the case in
the example of the SSE, that is why the LLN and the scaling limits
for a RW driven by these types of REs are serious mathematical
challenges. It is reasonable to think that if the observed unusual
regimes in Conjecture \ref{Scaling} will be rigorously proven, a
similar phenomenology can occur for other one dimensional ``slowly
mixing'' RE (e.g. other conservative particle systems) providing
that the RW is nearest-neighbor and has local drifts in both
directions. In principle, we might expect some analogous scenario
even in dimensions higher than one for static and dynamic REs
presenting again some strong correlation structure. The latter
observation is supported by the fact that a one-dimensional RW in
a dynamic RE can be seen as a two dimensional directed RW in a
static RE by interpreting the time as an extra spatial dimension
(see e.g. \cite{AvdHoRe11}).

We conclude with a table summarizing briefly some of the
one-dimensional rigorous and numerical results we presented so
far. Cells in gray mean that the corresponding result is
non-rigorous.

\vspace{1cm}

\newcommand{\cg}{\cellcolor[gray]{0.8}}
\begin{table}[h]
\centering
\begin{tabular}{c|c|c| rrrr}
\hline\hline
 &Static  &\multicolumn{2}{c}{Dynamic} \\ [0.5ex]
& &SSE (slowly mixing) & \;\;\;\;\;ISF (fast mixing)\\
\hline
Recurrence & $\rho=1/2$ & \cg $\rho=1/2$& $\rho=1/2$\;\;\;\;\;\;\\
\hline
LLN & & \cg $v>0$, & &&\\[-.5mm]
with & $v\geq0$ & \cg & $v>0$\;\;\;\;\;\;\;\;\\[-.5mm]
$\rho\neq1/2$ &  & \cg unclear for small $\gamma$& \\
\hline
& & \cg $\gamma$ big: diffusive, & \\[-.5mm]
Scaling & anomalous &  \cg & diffusive \;\;\;\;\;\;\\[-.5mm]
& & \cg $\gamma$ small: anomalous & \\
\hline
LDP rate $t$  & flat & flat & unique \;\;\;\;\;\;\;\\
for $v>0$& piece& piece&zero\;\;\;\;\;\;\;\;\;\; \\[1ex]
\hline\hline
\end{tabular}
\label{tab:hresult}
\end{table}

\vspace{1cm} {\bf{Acknowledgements:}} We would like to thank Frank
den Hollander, Renato dos Santos, Vladas Sidoravicius and Florian
V\"ollering for fruitful discussions during the preparation of
this paper. The simulations were performed on the Schr\"{o}dinger
cluster of the University of Z\"urich, where we gratefully
acknowledge the support of C. Bolliger, A. Godknecht, and R. Graf,
and at the Institute of Mathematics, thanks to R. Ostertag and C.
Rose.


\end{document}